\documentclass[lang = american]{ems-icm-arxiv} 

\usepackage{mathtools}
\usepackage{bbm}

\newtheorem{theorem}{Theorem}
\newtheorem{definition}[theorem]{Definition}
\newtheorem{conjecture}[theorem]{Conjecture}

\renewcommand{\C}{\mathbb{C}}
\newcommand{\F}{\mathbb{F}}
\newcommand{\I}{\mathbb{I}}

\newcommand{\N}{\mathbb{N}}
\renewcommand{\P}{\mathbb{P}}
\newcommand{\R}{\mathbb{R}}
\newcommand{\T}{\mathbb{T}}
\newcommand{\Z}{\mathbb{Z}}
\newcommand{\one}{\mathbbm{1}}

\newcommand{\II}{\mathcal{I}}
\newcommand{\BC}{\mathcal{BC}}

\renewcommand{\AA}{\mathcal{A}}
\newcommand{\BB}{\mathcal{B}}
\newcommand{\FF}{\mathcal{F}}
\newcommand{\GG}{\mathcal{G}}
\newcommand{\PP}{\mathcal{P}}

\renewcommand{\H}{\mathbf{H}}
\newcommand{\GP}{\mathbf{GP}}
\newcommand{\MM}{\mathbf{M}}

\newcommand{\e}{\mathsf{e}}
\newcommand{\f}{\mathsf{f}}

\newcommand{\Gr}{\mathrm{Gr}}
\newcommand{\Trop}{\mathrm{Trop\,}}
\newcommand{\PosTrop}{\mathrm{Trop^+}}
\newcommand{\cone}{\mathrm{cone}}
\renewcommand{\deg}{\mathrm{deg}}
\newcommand{\espan}{\mathrm{span}}

\newcommand{\csm}{\mathrm{csm}}
\newcommand{\MW}{\mathrm{MW}}
\newcommand{\Hom}{\mathrm{Hom}}

\numberwithin{equation}{section}

\begin{document}

\volumetitle{Federico Ardila--Mantilla} 

\title{The geometry of geometries: matroid theory, old and new}
\titlemark{The geometry of geometries: matroid theory, old and new}

\emsauthor{1}{Federico Ardila--Mantilla}{F.~Ardila}

\emsaffil{1}{San Francisco State University, San Francisco, California, USA; Universidad de Los Andes, Bogot\'a, Colombia \email{federico@sfsu.edu}}


\begin{abstract}
The theory of \emph{matroids} or \emph{combinatorial geometries} 
originated in linear algebra and graph theory, and has deep connections with many other areas, including 
field theory,
matching theory,
submodular optimization,
Lie combinatorics,
and total positivity. 
Matroids capture the combinatorial essence that these different settings share.

In recent years, the (classical, polyhedral, algebraic, and tropical) geometric roots of the field have grown much deeper, bearing  new fruits. 
We survey some recent successes, stemming from three geometric models of a matroid: the matroid polytope, the Bergman fan, and the conormal fan.

This survey was prepared for the 2022 International Congress of Mathematicians.
\end{abstract}

\maketitle


\section{Introduction}

There are natural notions of \emph{independence} in linear algebra, graph theory, field theory, matching theory, routing theory, rigidity theory, model theory, and many other areas.
When one seeks to understand the pleasing similarities between these different contexts, one is led to the powerful theory of \emph{matroids} or \emph{combinatorial geometries}. These intriguing, multifaceted objects turn out to also play a fundamental role in Lie combinatorics, tropical geometry, total positivity, and other settings.

The geometric approach to matroid theory has recently led to the solution of long-standing questions, and to the discovery of deep, fascinating interactions between 
combinatorics, algebra, and geometry. This survey is a selection of some recent achievements of the theory. 

Section \ref{sec:matroids} reviews basic definitions and Section \ref{sec:invariants} discusses key invariants of a matroid, like the characteristic and Tutte polynomials, and the finite field method to compute them.
Section \ref{sec:geom1} focuses on the \emph{matroid polytope}, starting from its parallel origins in combinatorial optimization and in the geometry of the Grassmannian. We discuss its connections to root systems, generalized permutahedra, the theory of matroid subdivisions and valuative invariants, and the role played by Hopf algebraic methods in these developments.
Section \ref{sec:geom2} concerns the \emph{Bergman fan} of a matroid. We discuss its central role in tropical geometry, 
its Hodge-theoretic properties, their role in the log-concavity of matroid $f$-vectors, and the theory of Chern-Schwartz-MacPherson cycles of matroids.
Section \ref{sec:geom3} discusses the \emph{conormal fan} of a matroid, its Hodge-theoretic properties, their role in the log-concavity of matroid $h$-vectors, and a Lagrangian geometric interpretation of CSM cycles. It also discusses two polytopes that play an important role in this theory and have elegant combinatorial properties: 
the bipermutahedron and harmonic polytope.
Finally, Section \ref{sec:future} offers some closing remarks.

\section{Matroids}\label{sec:matroids}

Matroids were defined independently in the 1930s by Nakasawa \cite{Nakasawa} and Whitney \cite{Whitney}. We choose one of many equivalent definitions. A \emph{matroid} $M=(E, \II)$ consists of a finite set $E$ and a collection $\II$ of subsets of $E$, called the \emph{independent sets}, such that

 (I1) $\emptyset \in \II$

 (I2) If  $J \in \II$ and $I \subseteq J$ then $I \in \II$.

 (I3) If $I, J \in \II$ and $|I|<|J|$ then there exists $j \in J-I$ such that $I \cup j \in \II$.

\noindent 
A matroid where every set of size at most 2 is independent is called a \emph{simple matroid} or a \emph{combinatorial geometry}.

Thanks to (I2), it is enough to list the collection $\BB$ of maximal independent sets; these are called the \emph{bases} of $M$. By (I3), they have the same size $r=r(M)$, which we call the \emph{rank} of $M$. Our running example will be the matroid with
\begin{equation}\label{eq:example}
E = abcde, \qquad \BB = \{abc, abd, abe, acd, ace\},
\end{equation}
omitting brackets for easier readability. Throughout the paper we let $n=|E|$ and $r=r(M)$.

The two most important motivating examples are graphical and linear matroids. Figure \ref{fig:exmatroid} shows how \eqref{eq:example} is a member of both families.

\smallskip

\noindent \textbf{\textsf{Graphical matroids.}} Let $E$ be the set of edges of a graph $G$ and $\mathcal{I}$ be the collection of forests of $G$; that is, the subsets of $E$ containing no cycle. 
\smallskip

\noindent \textbf{\textsf{Linear or realizable matroids.}} Let $\F$ be a field. \\
a) \textsf{(Vector configurations)}
Let $E$ be a finite set of vectors in a vector space over $\F$, and let $\II$ be the collection of linearly independent subsets of $E$. \\ 
b) \textsf{(Subspaces)} Let $V = \F^E$ be a vector space and $U \subseteq V$ be a subspace.
Let $\II$ be the collection of subsets $I \subseteq E$  such that $U$ intersects the coordinate subspace $V_I = \{\mathbf{v} \in V \, : \, v_i = 0 \text{ for } i \in I\}$ transversally; that is, $\dim (U \cap V_I) = \dim U - |I|$. 

\smallskip

The latter two constructions are equivalent: for a matrix $A$, the matroid of the set of columns of $A$ equals the matroid of the rowspace of $A$.

Matroids arise naturally in many important settings; \emph{e.g}, the study of algebraic dependences in a field extension, the combinatorics of root systems of semisimple Lie algebras, the perfect matchings in a bipartite graph, the non-intersecting paths in a directed graph, and total positivity of matrices, to name a few \cite{Oxley, Ardilasurvey}. 
Many of the matroids  in natural applications are linear, but   most matroids are not realizable over any field \cite{Nelson}, and including them leads to a much more powerful and robust theory of matroids and their geometry.

 \begin{center}
  \includegraphics[scale=1.4]{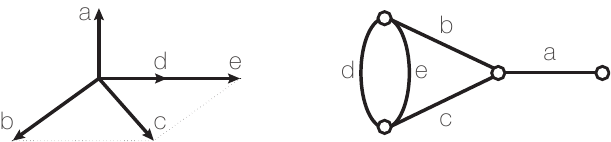}
  \captionof{figure}{\label{fig:exmatroid}
The matroid 
\eqref{eq:example}
with bases $\BB = \{abc, abd, abe, acd, ace\}$ is linear and graphical.}
 \end{center}

There are several natural operations on matroids. For $S \subseteq E$, the \emph{restriction} $M|S$ (or \emph{deletion} $M\backslash (E-S)$) and the \emph{contraction} $M/S$ are matroids on the ground sets $S$ and $E-S$, respectively, with independent sets
\begin{eqnarray*}
\II|S & \,\, \coloneq \,\,  & \{I \subseteq S \, : \, I \in \II\} \\
\II/S &\,\, \coloneq \,\, & \{I \subseteq E-S \, : \, I \cup I_S  \in \II\} 
\end{eqnarray*}
for any maximal independent subset $I_S$ of $S$; the latter is independent of the choice of $I_S$. When $M$ is a linear matroid in a vector space $V$, $M|S$ and $M/S$ are the linear matroids on $S$ and $E-S$ that $M$ determines on the vector spaces $\espan(S)$ and $V/\espan(S)$, respectively. If $S \subset T$ we write:
\[
M[S,T] \coloneq (M|T)/S.
\]

The \emph{direct sum} $M_1 \oplus M_2$ of $M_1 = (E_1, \II_1)$ and $M_2 = (E_2, \II_2)$ with $E_1 \cap E_2 = \emptyset$ is the matroid on  $E_1 \cup E_2$ with independent sets $\II_1 \oplus \II_2 = \{I_1 \cup I_2 \, : \, I_1 \in \II_1, I_2 \in \II_2\}$.
Every matroid decomposes uniquely as a direct sum of its \emph{connected components}.

Finally, the matroid $M^\perp$ \emph{dual} or \emph{orthogonal} to $M$ is the matroid on $E$ with bases
\[
\BB^\perp \coloneq \{E-B \, : \, B \in \BB\}.
\]
Remarkably, this simple notion simultaneously generalizes orthogonal complements and dual graphs.
If $M$ is the matroid for the columns of a matrix whose rowspan is $U \subseteq V$, then $M^\perp$ is the matroid for the columns of any matrix whose rowspan is $U^\perp$. If $M$ is the matroid for a planar graph $G$, drawn on the plane without edge intersections, then $M^\perp$ is the matroid for the dual graph $G^\perp$, whose vertices and edges correspond to the faces and edges of $G$, respectively, as shown in Figure \ref{fig:dualgraph}.

 \begin{center}
  \includegraphics[scale=.6]{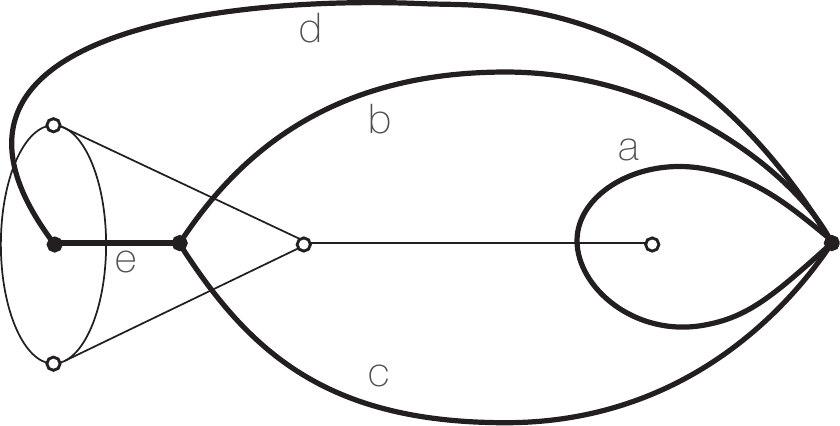}
  \captionof{figure}{ \label{fig:dualgraph}
  The dual matroid $\BB^\perp=\{bd, be, cd, ce, de\}$ is realized by the graph dual to Figure 1b.}
 \end{center}

An element $a$ is a \emph{loop} of $M$ if $\{a\}$ is dependent. A \emph{coloop} of $M$ is a loop of $M^\perp$.

\section{Enumerative invariants}\label{sec:invariants}

Two matroids $M_1 = (E_1, \II_1)$ and $M_2 = (E_2, \II_2)$  are \emph{isomorphic} if there is a \emph{relabeling} bijection $\phi:E_1 \rightarrow E_2$ that maps $\II_1$ to $\II_2$. 
A \emph{matroid invariant} is a function $f$ on matroids such that  $f(M_1) = f(M_2)$ whenever $M_1$ and $M_2$ are isomorphic. In 1964, Rota introduced the first foundational example \cite{RotaMobius}. 
\smallskip

\noindent \textbf{\textsf{The characteristic polynomial.}}
We define the \emph{rank function} $r: 2^E \rightarrow \Z$ of a matroid $M$ by
\[
r(A) \coloneq \textrm{ largest size of an independent subset of } A,
\]
for $A \subseteq E$. 
When $M$ is the matroid of a vector configuration $A$, $r(A) = \dim \espan(A)$.
The \emph{characteristic polynomial} of a loopless matroid $M$ is
\[
\chi_M(q) \coloneq \sum_{A \subseteq E} (-1)^{|A|} q^{r - r(A)}. 
\]

For Example \eqref{eq:example} we have $\chi_M(q) = q^3-4q^2+5q-2$. 
The characteristic polynomial of a matroid is one of its most fundamental invariants. For graphical and linear matroids, it has the following interpretations \cite{CrapoRota, OrlikSolomon, Zaslavsky}. 

\smallskip

\noindent  \textsf{1. Graphs.} If $M$ is the matroid of a connected graph $G$, then $q\, \chi_M(q)$ is the \emph{chromatic polynomial} of $G$; it counts the \emph{proper colorings} of the vertices of $G$ with $q$ given colors, where no two neighboring vertices have the same color.

\smallskip

\noindent \textsf{2. Hyperplane arrangements.}  Suppose $M$ is the matroid of  a set of non-zero vectors $v_1, \ldots, v_n$ spanning $\F^d$. Consider the arrangement $\AA$ of hyperplanes 
$H_i = \{x \in \F^d \, : \, v_i \cdot x = 0\}$ for $1 \leq i \leq n$, 
and its complement $V(\AA) = \F^d - (H_1 \cup \cdots \cup H_n)$. Depending on the underlying field, $\chi_M(q)$ stores different information about $V(\AA)$:
 
\noindent 
(a) ($\F = \C$) The Betti numbers of the complement $V(\AA)$ are the coefficients of $(-q)^d \chi_M(-1/q).$

\noindent 
(b) ($\F = \R$) The complement $V(\AA)$ consists of $(-1)^d\chi_M(-1)$ regions.

\noindent 
(c) ($\F = \F_q$) The complement  $V(\AA)$ consists of $\chi_M(q)$ points.
For a significantly stronger result on the $\ell$-adic \'etale cohomology of the arrangement, see \cite{BjornerEkedahl}.

\smallskip

Two related invariants that arise in several contexts are the \emph{M\"obius} and \emph{beta invariants} $\mu(M) = \chi_M(0)$ and $\beta(M) = (-1)^{r-1}\chi_M'(1)$, where $\chi'_M(x)$ is the derivative of $\chi_M(x)$. If $|E| \geq 2$ then  $\beta(M) = \beta(M^\perp)$.

\smallskip

\noindent \textbf{\textsf{The independence and broken circuit complex and their $f$- and $h$-vectors.}} Let $<$ be a linear order on $E$. 
A \emph{circuit} is a minimal dependent set. 
A \emph{broken circuit} is a set of the form $C - \{\min C\}$ where $C$ is a circuit. An \emph{nbc set} is a subset of $E$ not containing a broken circuit. 

We consider two simplicial complexes associated to a matroid $M$: the \emph{independence complex} and \emph{broken circuit complex}:
\[
\mathcal{I}(M) \coloneq \{\text{independent sets of } M\}, \qquad 
\mathcal{BC}_<(M) \coloneq \{\text{nbc sets of } M\}. 
\]

The \emph{$f$-vector} of a simplicial complex $\Delta$ of dimension $d-1$ counts the number $f_k(\Delta)$ of faces of $\Delta$ of 
size $k$ for each $0 \leq k \leq d$. The \emph{$h$-vector} of $M$ stores this information more compactly; it is given by
\[
\sum_{k=0}^{d} f_{k}(q-1)^{d-k} = 
\sum_{k=0}^{d} h_k q^{d-k}.
\]
The simplicial complexes
$\II(M)$ and $\BC_<(M)$
are $(r-1)$-dimensional. For Example \eqref{eq:example}, 
 \begin{eqnarray*}
f(\II(M)) = (1,5,9,5), & \qquad & f(\BC_<(M)) = (1,4,5,2), \\
h(\II(M)) = (1,2,2,0), &\qquad & h(\BC_<(M)) = (1,1,0,0).
\end{eqnarray*}
Topologically, $\II(M)$ and $\BC_<(M)$ are wedges of $\mu(M^\perp)$ and $\beta(M)$ spheres of dimension $r-1$, respectively \cite{Bjornermatroids}.
Up to alternating signs, the coefficients of $\chi_M(q)$ and $\chi_M(q+1)$ give the $f$-vector and $h$-vector of $\BC_<(M)$. 
\smallskip

\noindent \textbf{\textsf{The Tutte polynomial.}}
The invariant that appears most often in geometric, algebraic, and enumerative questions related to matroids is the \emph{Tutte polynomial} \cite{Tuttedichromatic}:
\begin{equation}\label{thm:Tutteformula}
T_M(x,y) \coloneq \sum_{A \subseteq M} (x-1)^{r-r(A)} \, (y-1)^{|A|-r(A)}.
\end{equation}
For Example \eqref{eq:example} we have $T_M(x,y) = x^3+x^2y+x^2+xy^2+xy$.

The ubiquity of this polynomial is explained by the following universality property. If a function $f:\textrm{Matroids} \rightarrow \F$ satisfies a \emph{deletion-contraction} of the following form:
\begin{equation} \label{eq:T-G}
f(M) =
\begin{cases} 
a f(M \backslash e) +  b f(M / e)  & \textrm{ if $e$ is neither a loop nor a coloop} \\
f(M \backslash e) \, f(L)  &  \textrm{ if $e=L$ is a loop}\\
f(M / e) \,  f(C)  & \textrm{ if $e=C$ is a coloop}
\end{cases}
\end{equation}
then it is an evaluation of the Tutte polynomial; namely, $f(M) = a^{n-r} \, b^r  \, T_M\left(\frac{f(C)}{b}, \frac{f(L)}{a}\right)$.
The Tutte polynomial also behaves very nicely under duality; we have $T_{M^\perp}(x,y) = T_M(y,x)$.

Many natural enumerative, algebraic, geometric, and topological quantities in numerous settings satisfy deletion-contraction recurrences, and hence are given by the Tutte polynomial; see \cite[Section 7.7.]{Ardilasurvey} and \cite{ArdilaTuttesurvey} for examples.
In particular, 
\[
\chi_M(q) = (-1)^r T_M(1-q,0),  \quad \mu(M) = T_M(1,0), \quad \beta(M) = [x^1y^0]T_M(x,y).
\]
where the last equality holds for $|E| \geq 2$, and
\begin{eqnarray*}
\sum_i f_{i}(\II(M)) x^{r-i} = T_M(1+x,1), & \qquad & \sum_i f_{i}(\BC_<(M)) x^{r-i} = T_M(1+x,0) \\
\sum_i h_i(\II(M)) x^{r-i} = T_M(x,1), &\qquad & \sum_i h_i(\BC_<(M)) x^{r-i} = T_M(x,0)
\end{eqnarray*}
In particular, though $\BC_<(M)$ depends delicately on the order $<$, its $f$- and $h$-vector do not.

\subsection{Computing Tutte polynomials: the finite field method}

Given how many quantities are given by the Tutte polynomial, it should not be a surprise that computing Tutte polynomials is extremely difficult (\#P-complete \cite{Welshcomplexity}) for general matroids. 
Nevertheless, we introduced a \emph{finite field method} \cite{ArdilaThesis, ArdilaTutte}, building on \cite{Athanasiadis, CrapoRota}, that has been effective for computing $T_M(x,y)$ in some special cases of interest. 
This method is similar in spirit to Weil's philosophy \cite{Weil} of learning about a complex projective variety $X(\mathbb{C})$ defined by  integer polynomials from its reductions to $X(\mathbb{F}_q)$ to various finite fields $\mathbb{F}_q$; in that setting we lose access to the complex  geometry but we gain the ability to count.

To compute the Tutte polynomial of $M$ with this method, we need a linear realization of $M$ as a set of vectors in $\mathbb{Q}^d$; most examples of interest have one.
For any power $q$ of a large enough prime number, this gives a linear realization of $M$ in $\F_q^d$; let ${\mathcal{A}_q}$ be the corresponding hyperplane arrangement. One then needs to count the points  in ${{\mathbb{F}}^d_q}$ according to the number of hyperplanes of $\mathcal{A}_q$ that they lie on. If one is able to do this for enough values of $q$, one can obtain the Tutte polynomial $T_M(x,y)$ from that enumeration. 

\begin{theorem} (Finite Field Method, \cite{ArdilaTutte}) Let $M$ be a matroid of rank $r$ realized  by a set of vectors in $\mathbb{Q}^d$. Let $q$ be a power of a large enough prime and $\mathcal{A}_q$ be the induced hyperplane arrangement in $\mathbb{F}_q^d$. Then
\[
\sum_{p \in {{\mathbb{F}}^d_q}} t^{h(p)} = q^{d-r}(t-1)^r T_M\left( 1+ \frac{q}{t-1}, t \right),
\]
where $h(p)$ is the number of hyperplanes of $\mathcal{A}_q$ containing $p$.
\end{theorem}

An equivalent result in the context of coding theory was obtained earlier by Greene \cite{Greeneweight}.
This finite field method is successful for root systems, arguably the most important vector configurations. 
The Tutte polynomial of the four families of \emph{classical root systems}
\[
A^+_{n-1} = \{\e_i - \e_j,\, : \,  i < j\}, \quad 
B(C)^+_n = \{\e_i \pm \e_j : \,   i <  j\} \cup \{(2)\e_i\}, \quad
D^+_n = \{\e_i \pm \e_j : \,   i <  j\}
\quad \subset \R^n,
\]
where $\{\e_1, \ldots, \e_n\}$ is the standard basis of $\R^n$, 
are given by the coefficient of $Z^n$ in the series
\[
T_A = F(Z,Y)^X,\qquad 
T_{BC} = F(2Z,Y)^{(X-1)/2}F(YZ,Y^2), \qquad 
T_D = F(2Z,Y)^{(X-1)/2}F(Z,Y^2),
\]
where $F(\alpha, \beta) = \sum_{n \geq 0} {\alpha^n \, \beta^{n \choose 2}}/{n!}$ is the \emph{deformed exponential function} \cite{ArdilaTutte}. Formulas for the exceptional root systems and the complex reflection groups are given in  \cite{DeConciniProcesi.Tutte, Geldon} and \cite{Randriamaro}.

The characteristic polynomial is particularly elegant: if $\Phi^+(\mathfrak{g})$ is the set of positive roots of a semisimple Lie algebra $\mathfrak{g}$ and $e_1, \ldots, e_n$ are its \emph{exponents}, then
\[
\chi_{\Phi^+(\mathfrak{g})}(q) = (q-e_1) \cdots (q-e_n).
\]
This is one of several examples where a characteristic polynomial surprisingly factors into linear factors; \cite{Saganfactors} outlines three conceptual explanations for this phenomenon.

The \emph{arithmetic Tutte polynomial} $M_A(x,y)$ of a vector configuration $A$ in a lattice also keeps track of arithmetic properties of $A$. This polynomial is related to the lattice point enumeration of zonotopes \cite{Stanleyzonotope, D'AdderioMoci.Ehrhart},  to complements of toric arrangements \cite{EhrenborgReaddySlone, Moci.toric, DeConciniProcesi.toric}, and to Dahmen-Micchelli and De Concini-Procesi-Vergne modules \cite{DahmenMicchelli, DeConciniProcesiVergne}. 
There is also a finite field method for computing $M_A(x,y)$  \cite{BrandenMoci, ArdilaCastilloHenley} that can be used successfully for root systems, and describes the volume and Ehrhart theory of Coxeter permutahedra \cite{ArdilaCastilloHenley, ArdilaBeckMcWhirter, DeConciniProcesi.Tutte}.

\section{Geometric model 1: Matroid polytopes}\label{sec:geom1}

A crucial insight on the geometry of matroids came from two seemingly unrelated places: combinatorial optimization and algebraic geometry. From both points of view, it is natural to model a matroid in terms of the following polytope. 

\begin{definition}
\cite{Edmonds}
The \emph{matroid polytope} of a matroid $M$ on $E$ is 
\[
P_M \coloneq \textrm{conv}\{\e_B  :  B \textrm{ is a basis of } M\} \,\, \subset \,\, \R^E,
\]
where  $\{\e_i  :  i \in E\}$ is the standard basis of $\R^E$ and  $\e_B = \e_{b_1} + \cdots + \e_{b_r}$ for $B =  \{b_1, \ldots, b_r\}$.
\end{definition}

Figure \ref{fig:matroidpolytope} shows this polytope for Example \eqref{eq:example}, a 3-dimensional polytope in $\R^5$.

 \begin{center}
  \includegraphics[scale=.8]{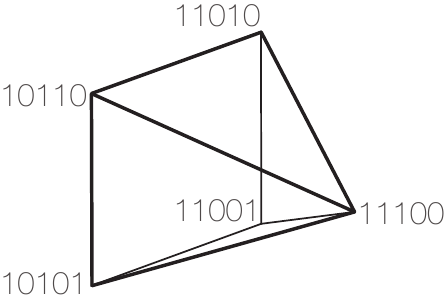}
  \captionof{figure}{ \label{fig:matroidpolytope}
The matroid polytope for the matroid \eqref{eq:example}. The vertices correspond to the bases.}
 \end{center}

\subsection{Algebraic geometry}

The intimate relation between matroids and the geometry of the Grassmannian is well-studied and mutually beneficial to both fields. Let us describe it briefly.

Instead of studying the $r$-dimensional subspaces of $\C^n$ one at a time, it is often useful to study them all at once. They can be conveniently organized into the \emph{Grassmannian} $\Gr(r,n)$; each point of $\Gr(r,n)$ represents an $r$-subspace of $\C^n$.

A choice of a coordinate system on $\C^n$ gives rise to the \emph{Pl\"ucker embedding} 
\begin{eqnarray*}
\Gr(r,n) &\,\,\, \xhookrightarrow{p} \,\,\, & \C\P^{{n \choose r}-1} \\
V & \,\,\, \longmapsto \,\,\, & \left(\det(A_B) \, : \, B \text{ is an $r$-subset of $[n]$} \right)
\end{eqnarray*}
defined as follows. For an $r$-subspace $V \subset \C^n$, choose an $r \times n$ matrix $A$ with $V = \mathrm{rowspan}(A)$. For each $r$-subset $B$ of $[n]$ let $p_B(V) \coloneq \det(A_B)$ be the determinant of the $r \times r$ submatrix $A_B$ of $A$ whose columns are given by the subset $B$.
Different choices of $A$ lead to the same \emph{Pl\"ucker  vector} $p(V)$ in projective space  $\C\P^{{n \choose r}-1}$.
The map $p$ provides a realization of the Grassmannian as a smooth projective variety.

Let $\C^* = \C \setminus \{0\}$. 
The torus $\T = (\C^*)^n/\C^*$ acts on the Grassmannian $\Gr(r,n)$ by stretching the $n$ coordinate axes of $\C^n$, modulo simultaneous stretching. Symplectic geometry then gives a \emph{moment map} $\mu: \Gr(r,n) \rightarrow \R^n$, which in this setting is given by
\[
\mu(V)  =\frac
{ \sum_{B}  |\det(A_B)|^2 \e_B}
{ \sum_{B}  |\det(A_B)|^2},
\]
where $B$ ranges over the $r$-subsets of $[n]$.

Now consider the orbit $\T \cdot V$ of the $r$-subspace $V \in \Gr(r,n)$ as the torus $\T$ acts on it, and the toric variety $X_V \coloneq \overline{\T \cdot V}$.
Gelfand, Goresky, MacPherson, and Serganova \cite{GGMS} proved that the moment map takes this toric variety to the matroid polytope of $M(V)$:
\begin{equation}\label{eq:moment}
\mu(\overline{\T \cdot V}) = P_{M(V)}.
\end{equation}
Thus matroid polytopes arise naturally in this algebro-geometric setting as well.

As a sample application, the degree of the closure of a torus orbit in the Grassmannian 
$X_V = \overline{\T \cdot V} \subset \C\P^{{n \choose r} -1}$ is then given by the volume of the matroid polytope $P_{M(V)}$. Ardila, Benedetti, and Doker \cite{ArdilaBenedettiDoker} used this to find a purely combinatorial formula for it. 

The projective coordinate ring of the toric variety $X_V$ is isomorphic to the subalgebra of the polynomial ring $\mathbb{C}[t_e \, : \, e \in E]$ generated by the monomials $t_B = \prod_{b \in B} t_b$ for the bases $B$ of $M(V)$. White \cite{Sturmfelstoriceqs, White80} conjectured that its defining toric ideal is generated by quadratic binomials. For the current state of the art on this conjecture, see \cite{Lason}.

\subsection{A geometric characterization of matroids}

 In most contexts where polytopes arise, it is advantageous if their faces can be described combinatorially. The vertices and edges often play an especially important role. For example, in geometry, they control the GKM presentation of the equivariant cohomology and K-theory of the Grassmannian \cite{FinkSpeyer, KnutsonTao}. In optimization, they are crucial to various algorithms for linear programming.

Matroid polytopes can be described entirely by their vertices and edges, as shown in the following beautiful combinatorial characterization. 

\begin{theorem} \label{thm:GGMS} \cite{GGMS}
A collection $\BB$ of subsets of $E$ is the set of bases of a matroid if and only if every edge of the polytope 
\[
P_\BB\coloneq {\mathrm{conv }}\{\e_B \, : \, B \in
{\mathcal{B}}\}\,\, \subset \,\, \R^E
\]
 is a translate of $\e_i-\e_j$ for some $i, j$ in $E$. 
\end{theorem}

Therefore, one could \emph{define} a matroid to be a lattice subpolytope of the cube $[0,1]^n$ that only uses these vectors as edges. 
Even if one is led to this family of polytopes through the geometry of subspaces as in \eqref{eq:moment}, one finds that non-linear matroids are equally natural from the polytopal point of view.
\textbf{Matroid theory provides the correct level of generality.}

 \begin{center}
  \includegraphics[scale=0.9]{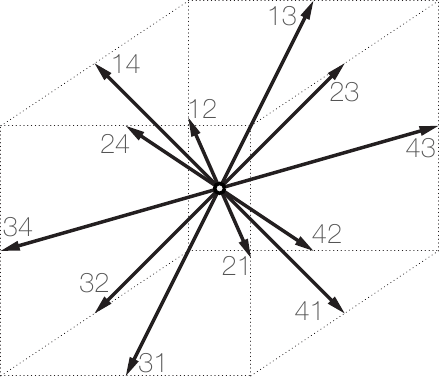}
  \captionof{figure}{ \label{fig:rootsystem} The root system $A_3 = \{\e_i-\e_j  :  1 \leq i, j \leq 4\},$ where  $\e_i-\e_j$ is denoted $ij$. Root systems play an essential role in matroid theory, as demonstrated by Theorem \ref{thm:GGMS}.}
 \end{center}

\noindent \textbf{\textsf{Positively oriented matroids}}.
\emph{Oriented matroids} are an abstraction of linear algebra over $\R$, abstracting linear dependence relations and their sign patterns. \emph{Positively oriented matroids} are those where every basis is positively oriented. Their matroid polytopes are precisely the polytopes whose edges are translates of $\e_i-\e_j$ and whose facet directions are $\overline{\e}_i - \overline{\e}_j$ for $i, j \in E$, where $\overline{\e}_i \coloneq \e_1 + \cdots + \e_i$. 
Ardila, Rinc\'on, and Williams \cite{ArdilaRinconWilliams2} used this characterization to prove da Silva's 1987 conjecture \cite{daSilva} that every positively oriented matroid is realizable.

\smallskip

\noindent \textbf{\textsf{Coxeter matroids}}.
Theorem \ref{thm:GGMS} shows that in matroid theory, a central role is played by one of the most important vector configurations in mathematics, the root system for the special linear group $SL_n$: 
\[
A_{n-1} = \{\e_i-\e_j \, : \, 1 \leq i, j \leq n\} \,\, \subset \,\, \R^n,
\]
 shown in Figure \ref{fig:rootsystem} for $n=4$. 
It is then natural to extend this construction to other semisimple Lie groups. The resulting theory of \emph{Coxeter matroids} \cite{Coxetermatroids}, introduced by Gelfand and Serganova, starts with a generalization of \eqref{eq:moment} and includes many other interesting results, but Borovik, Gelfand, and White's 2002 assessment still rings true today:
\begin{quote}
``the focal point of the theory: the relations between Coxeter matroids and the geometry of flag varieties [...] will need a few more years to settle in a definite form."  \cite{Coxetermatroids}
\end{quote}

\noindent
The enumerative combinatorics of Coxeter matroids, and its potential applications outside of matroid theory, are ripe for further exploration as well.

\subsection{Combinatorial optimization and generalized permutahedra} \label{sec:greedy}

Matroid theory also benefits from its close connection to submodular optimization, as discovered by Edmonds \cite{Edmonds} in 1970. This connection 
begins with a simple but fundamental observation: 
Matroids are precisely the simplicial complexes for which the greedy algorithm finds the facets of minimum weight. 

To be precise, for any \emph{weight function} $w: E \rightarrow \R$ on the elements of a matroid $M=(E,\II)$, let the weight of a basis $B$ be $w(B) = \sum_{b \in B} w(b)$. Then the bases of minimum weight are those obtained by applying the following greedy algorithm:

\begin{quote}
Start with $B=\emptyset$. Then, at each step, add to $B$ any element $e \notin B$ of minimum weight $w(e)$ such that $B \cup e \in {\mathcal{I}}$. Stop when $B$ is maximal in ${\mathcal{I}}$.
\end{quote}

\noindent Furthermore, a simplicial complex ${\mathcal{I}}$ constitutes the independent sets of a matroid if and only if this greedy algorithm works for 
any weight function $w$.

We may rewrite this greedy property as follows \cite{ArdilaKlivans, Coxetermatroids}.
Let $\mathcal{S}(w) \coloneq \{\emptyset = S_0 \subsetneq S_1 \subsetneq \cdots \subsetneq S_k \subsetneq S_{k+1} = E\}$ be the flag of subsets of $E$ such that $w(s_i)$ is constant for $s_i \in S_i - S_{i-1}$ and $w(S_1) < w(S_2 - S_1) < \cdots < w(S_k-S_{k-1}) < w(E - S_k)$. Then the $w$-minimum bases of $M$ are the bases of the \emph{$w$-minimum matroid} 
\begin{equation} \label{eq:M_w}
M_w \coloneq \bigoplus_{i=0}^k M[S_i, S_{i+1}].   
\end{equation}
It is useful to restate this geometrically as well.
The bases of the matroid $M$ are the vertices of the matroid polytope $P_M$. The $w$-minimum bases of $M$ are the vertices of the $w$-minimum face $(P_M)_w = \{x \in P_M \, : \, w(x) \leq w(y) \textrm{ for all } y \in P_M\}$ of $P_M$, and that face is itself the matroid polytope of $M_w$; that is, $(P_M)_w = P_{M_w}$.

Now consider the \emph{braid fan} $\mathcal{A}_E$ in $\R^E$ cut out by the hyperplanes $x_i = x_j$ for $i \neq j$ in $E$. This is the normal fan of the permutahedron $\Pi_E$, whose vertices are the $n!$ permutations of $[n]$.
The faces of $\mathcal{A}_E$ are in bijection with the flags of subsets of $E$: the open face $\sigma_\mathcal{S}$ consists of those  $w \in \R^E$ such that $\mathcal{S}(w) = \mathcal{S}$.
 Then 
\eqref{eq:M_w} shows that for any weight function $w$ in a fixed open face $\sigma_{\mathcal{S}}$, the matroid $M_w$ depends only on $\mathcal{S}$. This means that the braid fan $\mathcal{A}_E$ refines the normal fan of the matroid polytope $P_M$.

\smallskip

\noindent \textbf{\textsf{Generalized permutahedra and submodularity.}}
A \emph{generalized permutahedron}\footnote{Generalized permutahedra are the translates of the base polytopes of \emph{polymatroids}, of \cite{Edmonds}.} is a polytope $P$ in $\R^E$ satisfying the following three equivalent conditions \cite{Edmonds, Postnikovgenperm}:

$\bullet$ The braid fan $\mathcal{A}_E$ is a refinement of the normal fan of $P$.

$\bullet$ The edges of $P$ are parallel to roots $e_i-e_j$ for $i, j \in E$.

$\bullet$ $P = P(z) \coloneq \{x \in \R^E \, : \, \sum_{i=1}^n x_i = z(E), \sum_{i \in I} x_i \leq z(I)\}$ for a (unique) \emph{submodular function} $z$ on $E$: a function $z: 2^E \rightarrow \R$ with 
$z(A \cup B) + z(A \cap B) \leq z(A) + z(B)$ for  $A,B \subseteq E$.  \\
Allowing $z: 2^E \rightarrow \R \cup \{\infty\}$  we get \emph{extended generalized permutahedra}. 

The rank function of a matroid $M$ is one of the prototypical examples of a submodular function;  the corresponding generalized permutahedron is the matroid polytope $P_M$. 
This  explains why matroid theory informs and benefits greatly from the theory of submodular optimization.

\begin{figure}[!h]
\centering
\begin{center}
\includegraphics[scale=.3]{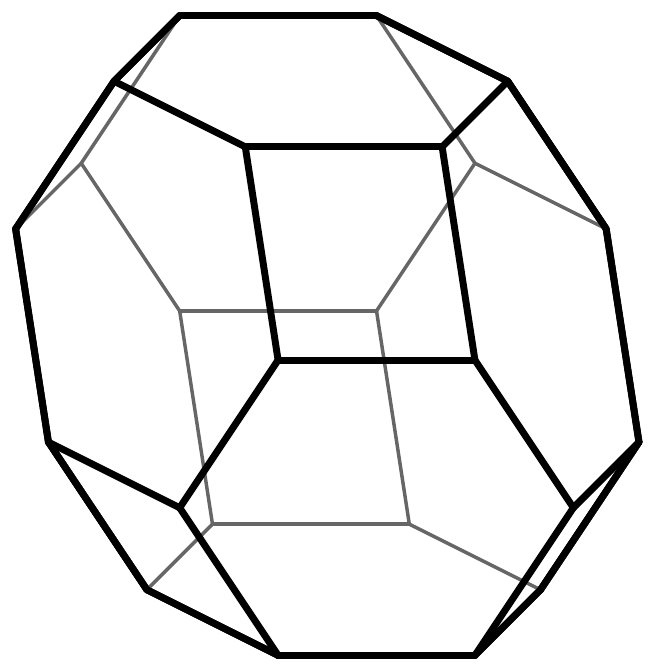} \qquad
\includegraphics[scale=.3]{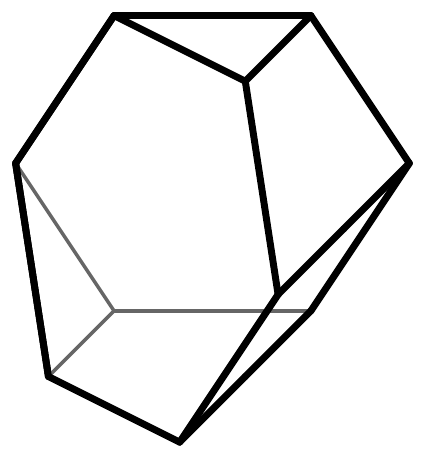}  \quad
\includegraphics[scale=.3]{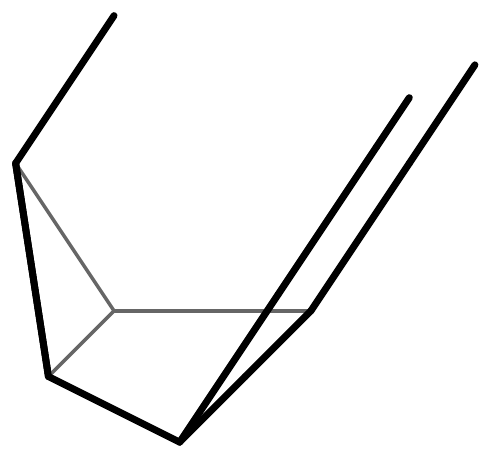} \quad
\includegraphics[scale=.3]{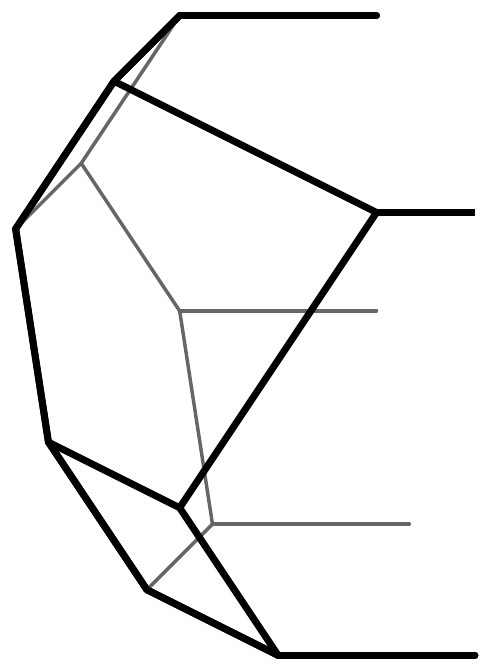} \quad
\includegraphics[scale=.3]{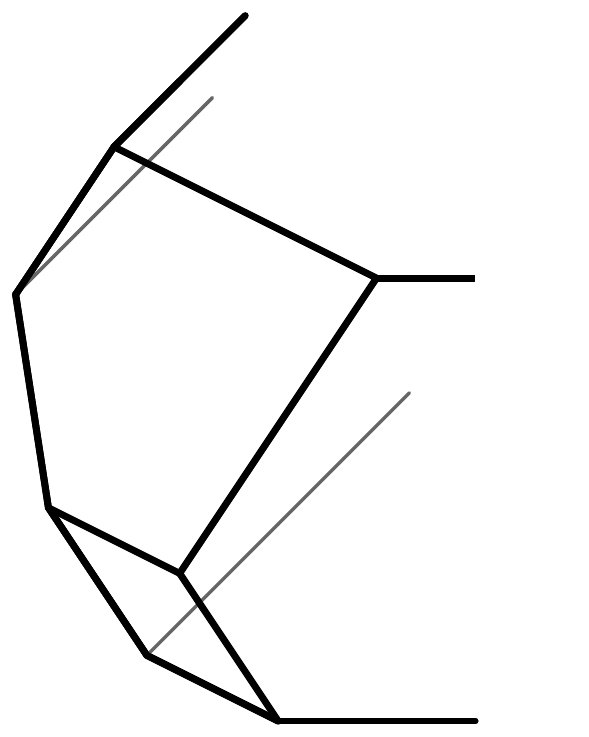}
\end{center}
\caption{The standard $3$-permutahedron and four other extended generalized permutahedra. \label{f:genperm}}
\end{figure}

Submodular functions arise in many contexts, partly because submodularity is equivalent to a natural \emph{diminishing returns property} that we now describe.
If $z$ measures some quantifiable benefit $z(A)$ associated to
each subset $A \subseteq E$, then the contraction $z/S(e) = z(S \cup e) - z(S)$ measures the \emph{marginal return} of adding $e$ to $S \notni e$. A function $z: 2^E \rightarrow \R$ is 
submodular if and only if $z/S(e) \geq z/T(e)$ for all $S \subseteq T \subseteq E-e$; that is, if the marginal return $z/S(e)$ diminishes as we add elements to $S \notni e$.

\subsection{Matroid subdivisions}

A \textbf{matroid subdivision} is a polyhedral subdivision $\mathcal{P}$ of a matroid polytope $P_M$ where every polytope $P \in \mathcal{P}$ is itself a matroid polytope. 
Equivalently, by Theorem \ref{thm:GGMS}, it is a subdivision of $P_M$ that only uses the edges of $P_M$. Let $\mathcal{P}^{int}$ be the set of interior faces of $\mathcal{P}$; these are the polytopes in $\mathcal{P}$ that are not on the boundary of $P_M$.
In the most important case, $M=U_{d,n}$ is the \emph{uniform matroid} where every $d$-tuple of $[n]$ is a basis, and  $P_M$ is the hypersimplex $\Delta(d,n)$. 

\begin{figure}[ht]
 \begin{center}
  \includegraphics[scale=.8]{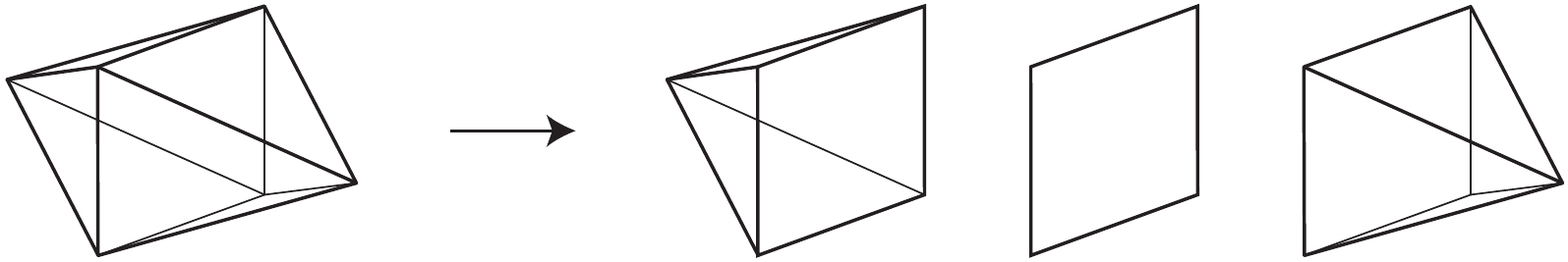}
  \caption{ \label{fig:matroidsubdiv}
The interior faces of a matroid subdivision of the uniform matroid $U_{2,4}$.}
 \end{center}
\end{figure}

Matroid subdivisions were first studied by Lafforgue in his 2003 work on surgery on Grassmannians \cite{Lafforgue}. These subdivisions also arose in algebraic geometry \cite{HackingKeelTevelev, Kapranov, Lafforgue}, in tropical geometry \cite{Speyer1}, and
in the theory of valuated matroids in optimization \cite{DressWenzel, Murota}.

Lafforgue gave an intriguing application of matroid subdivisions: 
if a matroid polytope $P_M$ has no nontrivial matroid subdivisions, then the matroid $M$ has (up to trivial transformations) only finitely many realizations over a fixed field ${\mathbb{F}}$. This is in stark contrast with Mn\"ev's Universality Theorem, which roughly states that every singularity type appears in the space of realizations of some oriented matroid $M$. This theorem was used by Vakil to construct several families of moduli spaces with arbitrarily bad singularities \cite{Vakil}.

The following conjecture of Speyer \cite{Speyer1} has led to many interesting developments. 

\begin{conjecture} \label{conj:fvector} If $M$ is a matroid on $[n]$ of rank $d$,
then a matroid subdivision of $M$ 
has at most $\frac{(n-c-1)!}{(d-c)!(n-d-c)!(c-1)!}$ interior faces of dimension $n-c$ for each $1 \leq c \leq \min\{d, n-d\}$. 
\end{conjecture}

\noindent 
For example, the subdivision of Figure \ref{fig:matroidsubdiv} has two interior faces of dimension three and one of dimension two, achieving equality in Conjecture \ref{conj:fvector}.

\subsection{Matroid valuations}\label{sec:valuations}

Matroid valuations are ways of measuring matroids that behave well under subdivision. Concretely, let $\mathsf{Mat}$ be the family of matroids and $A$ be an abelian group. A function $f:\mathsf{Mat} \rightarrow A$ is a \emph{matroid valuation} if for any subdivision $\PP$ of a matroid polytope $P_M$ we have the inclusion-exclusion relation
\begin{equation}\label{eq:weakval}
f(M) = \sum_{P_N \in \, \PP^{int}} (-1)^{\dim P_M  - \dim P_N} f(N). 
\end{equation}

The volume, the number of lattice points, and the Ehrhart polynomial (given by $\mathrm{Ehr}_P(t) = |tP \cap \Z^d|$ for $t \in \N$) are natural ways of measuring a polytope, and an Euler characteristic computation shows that they are matroid valuations.
More interestingly, matroids can also be measured using seemingly unrelated combinatorial and algebro-geometric invariants that, unexpectedly, also satisfy \eqref{eq:weakval}. These valuations include, among many others:

$\bullet$ the Tutte polynomial of a matroid \cite{Speyer1, ArdilaFinkRincon},

$\bullet$ the Chern-Schwartz-MacPherson cycles of a matroid \cite{LRS},

$\bullet$ the Kazhdan-Lusztig polynomial of a matroid \cite{EPW, ArdilaSanchez},

$\bullet$ the motivic zeta function of a matroid \cite{JKU, ArdilaSanchez},

$\bullet$ the Speyer polynomial of a matroid \cite{Speyer2}.

$\bullet$ the volume polynomial of the Chow ring of a matroid \cite{Eurvolume}.

\noindent 
Ardila, Fink, and Rinc\'on \cite{ArdilaFinkRincon} gave a general geometric technique to construct many valuations of matroids, including valuative invariants.
Derksen and Fink \cite{DerksenFink} constructed the universal valuative invariant of matroids. In practice, it is often not clear why a conjectural valuation is a function of the Derksen-Fink invariant.

As a sample application of matroid valuations, we sketch Speyer's proof of Conjecture \ref{conj:fvector} for $c=1$. Since $T_M(x,y)$ is a valuation, so is the \emph{beta invariant} $\beta(M) = [x^1y^0]T_M(x,y)$. The deletion-contraction recursion gives $\beta(U_{d,n}) = {n-2 \choose d-1}$, and one can show that 
$\beta(N)=0$ if $N$ is not connected (or, equivalently, if $P_N$ is not full-dimensional) and $\beta(N) \geq 1$ otherwise.  
Thus, for a matroid subdivision $\PP$ of $U_{d,n}$,
\[
{n-2 \choose d-1} = \beta(U_{d,n}) = \sum_{\stackrel{P_N \in \,  \mathcal{P}}{P_N  \textrm{ facet}}}\beta(N) \geq (\textrm{number of facets of } \mathcal{P}).  
\]

Similarly, 
Speyer \cite{Speyer2} constructed a polynomial invariant $g_M(t)$ motivated by the K-theory of the Grassmannian, and he used it to prove Conjecture \ref{conj:fvector} for matroid subdivisions whose matroids are realizable over a field of characteristic $0$. 
His proof relies on the nonnegativity of $g_M(t)$, which is only known for matroids realizable in characteristic $0$, for which the coefficients of this polynomial have a geometric interpretation.
The nonnegativity of $g_M(t)$ for all $M$, which would prove Conjecture \ref{conj:fvector} in full generality, remains open.
 
\smallskip

The numerous examples of this section raise a natural question. 
In the following section we offer a possible answer coming from Hopf algebras.

\smallskip

\noindent 
\textbf{Question.} Why are many natural functions of matroids also matroid valuations? How might we find others?

\smallskip

\subsection{Hopf algebras and valuations}

In 1978, Joni and Rota \cite{JoniRota} showed that many combinatorial families have natural \emph{merging} and \emph{breaking} operations that give them the structure of a Hopf algebra, with many useful consequences. 
For the Hopf algebra of matroids, a pleasant surprise was uncovered recently: the geometric point of view plays a central role, and connects naturally with the theory of matroid valuations.

\smallskip

\noindent \textbf{\textsf{The Hopf algebra of matroids.}}
Joni--Rota  \cite{JoniRota} and Schmitt \cite{Schmitt} defined the \emph{Hopf algebra of matroids} $\MM$ as the span of the set of matroids modulo isomorphism, with the product $\cdot: \MM  \otimes \MM \rightarrow \MM$ and coproduct $\Delta: \MM \rightarrow \MM \otimes \MM$ given by:
\[
M \cdot N  \coloneq M \oplus N, 
\qquad  \qquad
\Delta(M)  \coloneq \sum_{S \subseteq E} (M|S) \otimes (M/S). 
\]
A Hopf algebra has an \emph{antipode map} $S:\MM \rightarrow \MM$, which is the Hopf-theoretic analog of the inverse map $g \mapsto g^{-1}$ in groups. 
Takeuchi \cite{Takeuchi} gave a general formula for the antipode of any connected, graded Hopf algebra; it is an alternating sum with a superexponential number of terms, that is generally not tractable.
Thus a central problem for a Hopf algebra of interest  $H$ is to use the structure of $H$ to find an explicit, cancellation-free formula for $S$. 

The optimal formula for the antipode of matroids was  discovered by Aguiar and Ardila \cite{AguiarArdila}. 
The key new insight is that, although they arose in optimization and geometry,
\textbf{matroid polytopes are also fundamental in the Hopf algebraic structure of matroids.}

\begin{theorem} \label{thm:antipode} \cite{AguiarArdila}
The antipode of the Hopf algebra of matroids $\MM$ is 
\[
S(M) = \sum_{P_N \textrm{ face of } P_M} (-1)^{c(N)} N
\]
for any matroid $M$, where $c(N)$ denotes the number of connected components of $N$.
\end{theorem}

\smallskip

\noindent \textbf{\textsf{The incidence Hopf algebra of matroids.}} Theorem \ref{thm:antipode} makes it very tempting to replace each matroid $M$ with the indicator function $\one_{M} : \R^E \rightarrow \R$ given by $\one_{M}(p) = 1$ if $p$ is in the matroid polytope $P_M$ and  
$\one_{M}(p) = 0$ otherwise. This would give the formula
\[
S(P_M) = (-1)^{c(M)} \text{int} (P_M),
\]
suggesting connections with the \emph{Euler map} of McMullen's polytope algebra \cite{McMullen} and with \emph{Ehrhart reciprocity} for lattice polytopes \cite{Ehrhart}.

Ardila and Sanchez made this precise, constructing the \emph{incidence Hopf algebra of matroids} $\I(\MM)$. Its component of degree $n$ is spanned by the indicator functions of matroid polytopes on $[n]$. We have
\[
\I(\MM) \cong \MM/\mathrm{ie}(\MM) 
\text{ \quad for \quad } 
\mathrm{ie}(\MM) \coloneq \text{span} \{P_M - 
\sum_{P_N \in \, \PP^{int}} (-1)^{\text{codim} P_N}  P_N \, : \, \PP \text{ subdivides }P_M\}.
\]
The subspace $\mathrm{ie}(\MM)$ is a Hopf ideal of $\MM$, so the quotient $\I(\MM)$ is indeed a Hopf algebra. 

Matroid valuations are precisely the functions on matroids that descend to the vector space $\I(\MM)$. The Hopf algebraic structure on $\I(\MM)$ then provides a straightforward, unifying framework to discover and prove many known and new matroid valuations, including all the ones discussed in Section \ref{sec:valuations}. This offers one possible answer to Question 1.

\smallskip

\noindent \textbf{\textsf{Generalized permutahedra and universality.}} The constructions and theorems of this section hold more generally for the Hopf algebras
 $\GP^{(+)}$ of (extended) generalized permutahedra and their indicator functions. The Hopf algebras of symmetric functions, Fa\'a di Bruno, matroids, graphs, posets, and many others can be realized as subalgebras of $\GP^+$, with many useful consequences. 
 The \emph{character theory} of Hopf algebras and Theorem \ref{thm:antipode} give a unified explanation of various \emph{combinatorial reciprocity theorems}: instances where the same polynomial $p(x)$ gives the number $p(n)$ of $A$-structures of size $n$ and the number $|p(-n)|$ of $B$-structures of size $n$ for two \emph{reciprocal} combinatorial families $A$ and $B$.
 They also explain  why the coefficients of the multiplicative and compositional inverses of power series are given by the face structure of the permutahedron and associahedron, respectively \cite{AguiarArdila}.
 
This raises another natural question:

\smallskip

\noindent 
\textbf{Question 2.} Why are many Hopf algebras in combinatorics related to generalized permutahedra and their characters?

As a partial answer to this question, we offer two universality results.

\medskip

To define $\GP$, the key fact is that for any generalized permutahedron $P \subset \R^E$ and any subset $S \subseteq E$, the $\e_S$-maximal face of $P$ decomposes as  $P_S = (P|S) \times (P/S)$ for generalized permutahedra $P|S \subset \R^S$ and $P/S \subset \R^{E-S}$. Define a product and coproduct by 
\begin{equation} \label{eq:Hopf}
P \cdot Q \coloneq P \times Q, 
\qquad \qquad
\Delta(P) \coloneq \sum_{S \subseteq E} (P|S) \otimes (P/S).
\end{equation}
Aguiar and Ardila \cite{AguiarArdila} proved that generalized permutahedra are the maximal family of polytopes for which \eqref{eq:Hopf} defines a Hopf algebra. The antipode is analogous to Theorem \ref{thm:antipode} \cite{AguiarArdila}.

In a different, and more general direction, we have the following universality theorem. 
A \emph{generalized polynomial character} on a Hopf algebra $\H$ is a multiplicative function from $\H$ to the ring of generalized polynomials, which can have any real numbers as exponents.
For example, the \emph{canonical character} on $\GP^+$ is $\beta(P) = t^{r(P)}$ where $P$ lies on the hyperplane $\sum_{i \in E} x_i = r(P)$.
Ardila and Sanchez \cite{ArdilaSanchez} proved that the indicator Hopf algebra 
$(\I(\GP^+),\beta)$ is the terminal Hopf algebra with a generalized polynomial character.
\textbf{Any Hopf algebra with a generalized polynomial character factors through $\I(\GP^+)$.}
This partially explains the ubiquity of these polytopes in combinatorial Hopf algebras.

These results are closely related to Derksen and Fink's universal valuative invariant for generalized permutahedra \cite{DerksenFink}. Generalizing to the setting of finite root systems, Ardila, Castillo, Eur, and Postnikov described \emph{generalized Coxeter permutahedra} \cite{ArdilaCastilloEurPostnikov} and Eur, Sanchez, and Supina computed their universal valuation \cite{EurSanchezSupina}. It would be interesting to construct a Coxeter-Hopf-theoretic framework where generalized Coxeter permutahedra, Coxeter matroids, and other related objects fit naturally.

\section{Geometric model 2: Bergman fans}
\label{sec:geom2}

We now introduce a second geometric model of matroids, coming from tropical geometry. 
It relies on the \emph{flats} of $M$; these are the subsets $F \subseteq E$ such that 
$
r(F \cup e) > r(F) \textrm{ for all } e \notin F.
$
We say a flat $F$ is \emph{proper} if it does not have rank $0$ or $r$.
The \emph{lattice of flats} of $M$, denoted $L_M$, is the set of flats, partially ordered by inclusion.

 \begin{center}
    \includegraphics[height=1.5cm]{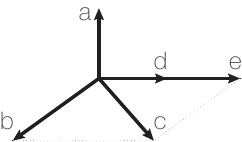} \qquad 
  \includegraphics[height=3cm]{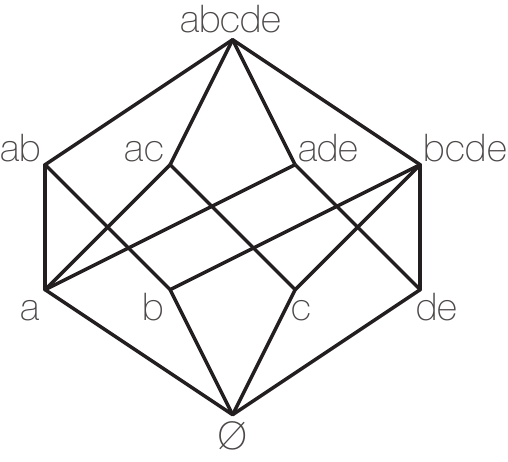} \qquad 
    \includegraphics[height=3.5cm]{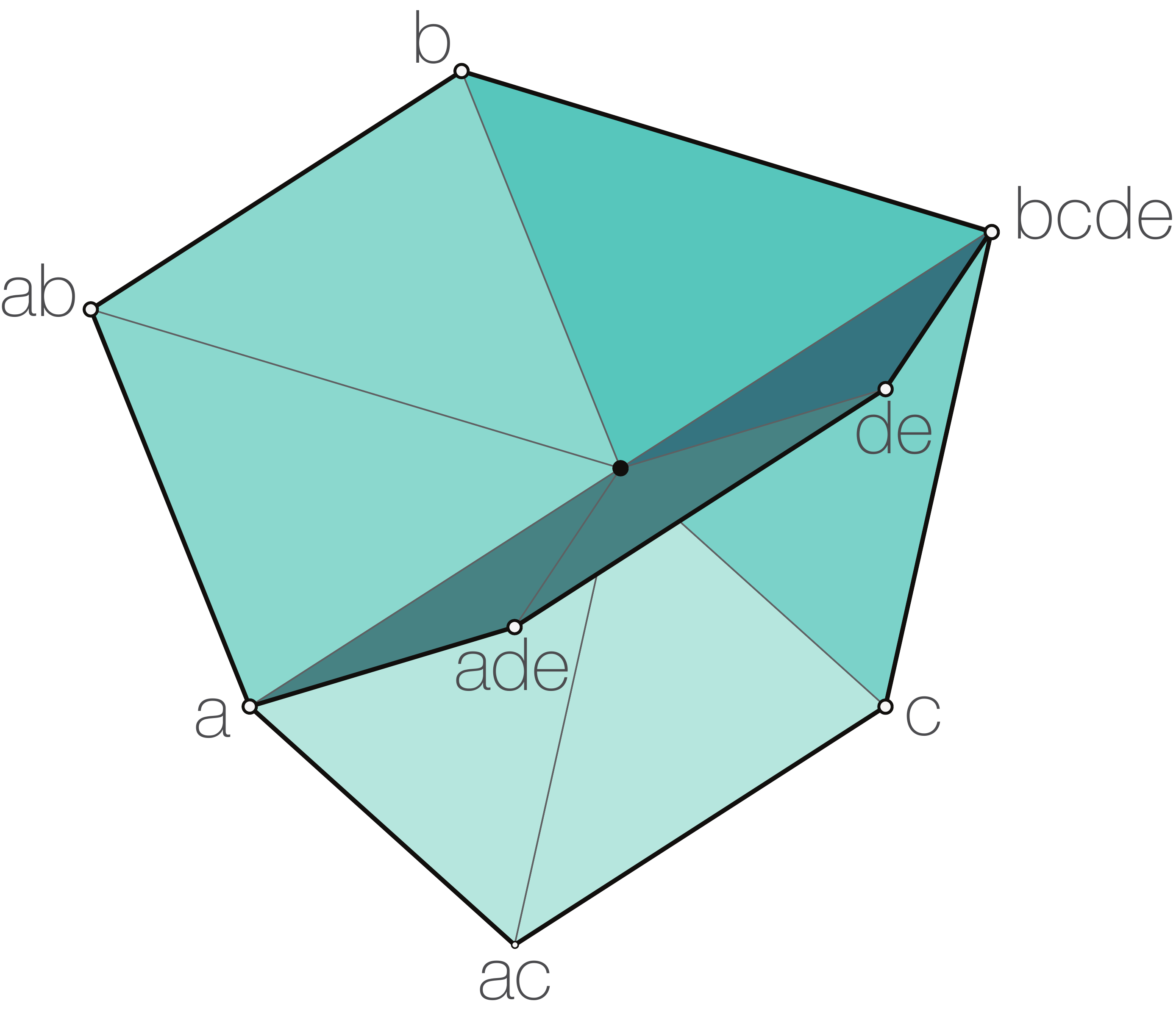}
  \captionof{figure}{  \label{fig:lattice}
Our sample matroid \eqref{eq:example}, its lattice of flats, and its Bergman fan, which is the cone over a wedge of $|\mu(M)| = 2$ spheres.}
 \end{center}

When $M$ is the matroid of a vector configuration $E$ in a vector space $V$,  the flats of $M$ correspond to the subspaces of $V$ spanned by subsets of $E$, as illustrated in Figure \ref{fig:lattice}. In this section we assume that the matroid $M$ has no loops.

Ardila and Klivans introduced the following polyhedral rendering of a matroid:

\begin{definition}\label{def:Bergman} \cite{ArdilaKlivans}  
The \emph{Bergman fan} or \emph{matroid fan} $\Sigma_M$ of a matroid $M$ on $E$ is the polyhedral fan in $\R^E/\left< \e_E \right>$ consisting of the cones
\[
\sigma_\FF \coloneq \cone \{\e_F \, : \, F \in \FF \} 
\]
for each flag 
$\FF = \{F_1 \subsetneq \cdots \subsetneq F_l\}$ of proper flats of $M$.
Here $\e_F \coloneq \e_{f_1} + \cdots + \e_{f_k}$ for $F= \{f_1, \ldots, f_k\}$. \\
\end{definition}

Let us discuss the tropical geometric origin of this fan, and some of its applications.

\subsection{Tropical geometry}   Tropical geometry is a powerful technique\footnote{whose name is questionable to this mathematician from the tropics}  designed to answer questions in algebraic geometry by translating them into polyhedral questions that can be approached combinatorially. 
In one of its manifestations, tropical geometry sends a complex algebraic variety $V \subset (\C^*)^n$ to its amoebas $\AA_t(V)$, whose limit as $t$ approaches $0$ is a piecewise linear space $\Trop V$
called the \emph{tropicalization} 
or \emph{logarithmic limit set} of $V$:
\[
\AA_t(V) \coloneq \{(\log_{\frac1t}(|z_1|, \ldots \log_{\frac1t}(|z_n|)) \, : \, (z_1, \ldots, z_n) \in V\}, \qquad 
\Trop V \coloneq \lim_{t \rightarrow 0} \mathcal{A}_t(V). 
\]
For an introduction and a more precise discussion, see \cite{Alessandrini, AllermannRau, MaclaganSturmfels, MikhalkinRau}.
An important early success of the theory was Mikhalkin's 2005 tropical computation \cite{Mikhalkin} of the \emph{Gromov-Witten invariants of $\C\P^2$}, which count the plane curves of degree $d$ and genus $g$ passing through $3d+1-g$ general points.  Since then, many new results in classical algebraic geometry have been obtained through tropical techniques.

The tropical approach requires two steps. Firstly, one needs to recognize what features of a geometric object $V$ can be recovered from its tropicalization $\Trop V$, which only captures part of the behavior of $V$ at infinity. Secondly, one needs to realize that $\Trop V$ may be simpler than $V$, but it is still usually very intricate. 

Additionally, to develop a robust theory, one is led to define \emph{tropical varieties} that are not necessarily tropicalizations of algebro-geometric objects, but are equally important tropically.
Understanding their structure is the source of very interesting combinatorial problems. 

An important development towards the algebraic foundations of tropical geometry was Sturmfels's description \cite{Sturmfelspolyeqs} of 
$\Trop V$ in terms of the Gr\"obner fan of its ideal $I(V)$:
\[
\Trop V = \{w \in \R^n \, : \, \textrm{the $w$-initial ideal of $I(V)$ has no monomials}\}.
\]
This led him to define the \emph{tropical variety of a matroid} $M$ on $E$ to be
\[
\Trop M \coloneq \{w \in \R^E \, : \, \textrm{the $w$-minimum matroid } M_w \textrm{ has no loops}\}.
\]
 If $M=M(V)$ is the matroid of a linear subspace $V \subset (\C^*)^n$ then $\Trop M = \Trop V$. If $M$ is not linear, $\Trop M$ is not the tropicalization of a variety.

Bergman \cite{Bergman} conjectured and  Bieri and Groves \cite{BieriGroves} showed that the tropicalization of an irreducible variety in $(\C^*)^n$  is pure and connected. Sturmfels \cite{Sturmfelsconjecture} conjectured that $\Trop M$ should have these same properties, even if $M$ is not a linear matroid.  
Ardila and Klivans \cite{ArdilaKlivans} first introduced the Bergman fan with the goal of settling this conjecture.

\begin{theorem} \label{thm:AK} 
The Bergman fan $\Sigma_M$ of a matroid $M$ is a triangulation of the tropical space $\Trop M$. Therefore $\Trop M$ is a cone over a wedge of $|\mu(M)|$ spheres of dimension $r-2$, where $\mu(M)$ is the M\"obius number of the matroid.
\end{theorem}

The first statement follows from the fact, explained in Section \ref{sec:greedy}, that the matroid $M_w$ only depends on the face of the braid fan $\AA_E$ containing $w$. Intersecting $\Trop M$ with the braid fan induces the triangulation $\Sigma_M$.  
The second statement then follows from Bj\"orner's result \cite{Bjornermatroids} that the order complex of the lattice of flats of $M$ is a wedge of spheres. 

\smallskip

\noindent \textbf{\textsf{Total positivity and oriented matroids}}.
Motivated by the theory of total positivity, Speyer and Williams \cite{SpeyerWilliams} introduced the positive part $\PosTrop V$ of the tropicalization of an affine variety $V$. Analogously, Ardila, Klivans, and Williams \cite{ArdilaKlivansWilliams} studied the \emph{positive part of the tropical variety of an acyclic oriented matroid} $M$, which is 
$\PosTrop M = 
 \{w \in \R^E \, : \, \textrm{the $w$-minimum matroid } M_w \textrm{ is acyclic}\}$.
They showed that 
the \emph{LasVergnas face lattice} of $M$ gives a triangulation of $\PosTrop M$; and hence \cite{OMs}  $\PosTrop M$  is a cone over a sphere.

Furthermore, 
Ardila, Reiner, and Williams \cite{ArdilaReinerWilliams}
 constructed 
$|\mu(M)|$ reorientations $M^\epsilon$ of $M$ that decompose $\Trop(M)$ explicitly as a wedge of the corresponding spheres $\PosTrop M^\epsilon$  \cite{ArdilaReinerWilliams}.
When $M$ is the matroid of a root system $\Delta$, each sphere $\PosTrop M^\epsilon$ is dual to the \emph{graph associahedron} of the Dynkin diagram of $\Delta$ \cite{ArdilaReinerWilliams}.

There are also very interesting connections between the positive part of the tropical Grassmannian, matroid theory, and cluster algebras; see \cite{WilliamsICM} for a survey.

\subsection{A tropical characterization of matroids}

A \emph{tropical fan} is a subset $X \subseteq \R^n$ that has the structure of a pure, integral polyhedral fan $\mathcal{X}$ with weights $w:\{\text{facets of } \mathcal{X}\} \rightarrow \N$ satisfying the \emph{balancing condition}:
\begin{equation} \label{eq:balancing}
\sum_{\textrm{facets } \sigma \supset \tau} w(\sigma) \, v_{\sigma/\tau} = 0 
\textrm{  mod } \espan(\tau) \qquad
\textrm{ for any face $\tau$ of codimension $1$,}
\end{equation}
where $v_{\sigma/\tau} \in \Z^n$ is the primitive generator of the ray $\sigma/\tau$ in $\Z^n / \Z \tau$. Examples include $\Trop V$ for any subvariety $V \subset (\C^*)^n$ and $\Trop M$ for any matroid $M$ on $[n]$.

In analogy with the classical setting, the \emph{degree} of a tropical fan $X$ is obtained by counting the intersection points of $X$ with $\Trop V$ for a generic linear subspace $V$ of codimension $\dim X$, 
with certain multiplicities. 
For precise definitions, see \cite{AllermannRau, MikhalkinRau}.\footnote{One sometimes allows bounded faces in a tropical variety; Fink works in this setting.}
Fink \cite{Fink} gave the following remarkable characterization: 

\begin{theorem} \label{thm:charBergman} \cite{Fink}
The tropical fans of degree 1 in $\R^E$ are precisely the tropical varieties of the matroids on $E$.
\end{theorem}

Thus Bergman fans of matroids can be thought of as the tropical analogs to linear subspaces. 
In fact, one could \emph{define} a matroid on $E$ to be a tropical fan of degree $1$ in $\R^E$. 
 Notice that, although the matroid fan $\Trop M$ only arises via tropicalization when $M$ is a linear matroid, one should really consider the matroid fans of non-realizable matroids as well; they are equally natural from the tropical point of view.
Again, \textbf{matroid theory  provides the correct level of generality}.

Theorems \ref{thm:AK} and \ref{thm:charBergman} explain two important roles that matroids play in tropical geometry. 
On the one hand, they offer a useful testing ground, providing hints for the kinds of general results that may be possible, and the sorts of difficulties that one should expect. 
On the other hand, they are fundamental building blocks; for instance, in analogy with the classical definition of a manifold, a \emph{tropical manifold}
is an abstract tropical variety that locally has the structure of a Bergman fan of a matroid \cite{MikhalkinRau}; Figure \ref{fig:manifold} shows an example. Clarifying the foundations of the theory of tropical manifolds is an important project; we expect it  will continue to shape and benefit from the further development of the geometry of matroids.

 \begin{center}
  \includegraphics[scale=.2]{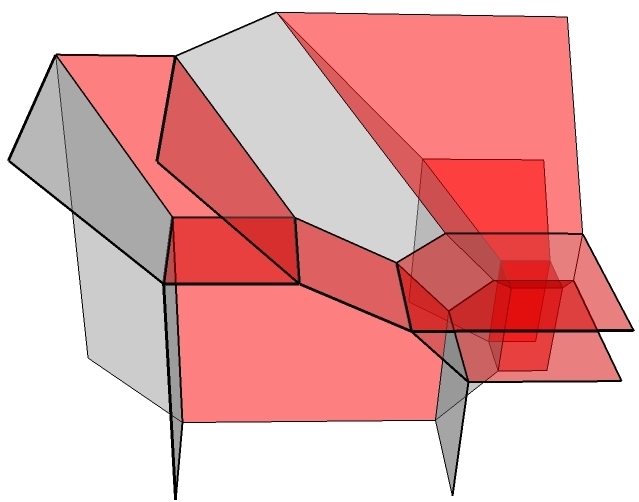}
  \captionof{figure}{ \label{fig:manifold} A tropical manifold is a tropical variety that is a matroid fan locally. (Picture: Johannes Rau)}
 \end{center}

Another interesting direction is the tropical geometry of Coxeter matroids and homogenous spaces; for some initial efforts, see \cite{BallaOlarte, BrandtEurZhang, RinconD}. As we will see in the next section, this could have interesting enumerative applications.

\subsection{The Chow ring, combinatorial Hodge theory, and log-concavity}

We say a sequence $a_0, a_1, \ldots, a_r$ of non-negative integers is: 

\vspace{.1cm}
\noindent 
$\bullet$  \emph{unimodal} if $a_0 \leq a_1 \leq \cdots \leq a_{m-1} \leq a_m \geq a_{m+1} \geq \cdots  \geq a_r$ for some $0 \leq m \leq r$, \\
\noindent 
$\bullet$ \emph{log-concave} if $a_{i-1}a_{i+1} \leq a_i^2$ for all $1 \leq i \leq r-1$, and \\
\noindent 
$\bullet$ \emph{flawless} if $a_i \leq a_{s-i}$ for all $1 \leq i \leq \frac{s}2$, where $s$ is the largest index with $a_s \neq 0$. 
\vspace{.1cm}

\noindent Many sequences in combinatorics have these properties, but proving it 
often requires a fundamentally new construction or connection to algebra or geometry, and gives rise to unforeseen structural results about the objects of interest. 

In 1970, Rota \cite{RotaICM} first raised such questions in the context of matroid theory, and suggested the Alexandrov-Fenchel inequalities in convex geometry as a possible approach. In the early 80s, Stanley \cite{StanleyICM, StanleyLie} systematically used the hard Lefschetz theorem and the representation theory of Lie algebras to prove similar combinatorial inequalities. In recent years, building on these techniques, a \emph{combinatorial Hodge theory of matroids} has led to the solution of several long-standing open problems in matroid theory.

The \emph{Chow ring} of the Bergman fan $\Sigma_M$ is
\[
A^*(\Sigma_M) \coloneq \R[x_F \, : \, F \textrm{ proper flat of } M]/(I_M + J_M),
\]
where
\[
I_M \coloneq \left< x_{F_1}x_{F_2} \, : \, F_1 \subsetneq F_2 \textrm{ and } F_1 \supsetneq F_2\right>,  \qquad 
J_M \coloneq \left<\sum_{F \ni i} x_F - \sum_{F \ni j} x_F \, : \, i, j \in E\right>. 
\]
Work of Brion \cite{Brion} implies that $A^*(\Sigma_M)$  is isomorphic to the Chow ring of the toric variety associated to $\Sigma_M$. Work of Billera \cite{Billera} implies that $A^*(\Sigma_M)$ is also isomorphic to the algebra of continuous piecewise polynomial functions on $\Sigma_M$, modulo the restrictions of global linear functions to $\Sigma_M$. When studying this ring, it is often useful to keep in mind both its algebraic presentation  and its interpretation in terms of piecewise polynomial functions.

When $M$ is linear over $\C$, Feichtner and Yuzvinsky \cite{FeichtnerYuzvinsky} proved that  $A^*(\Sigma_M)$ is the Chow ring of De Concini and Procesi's \emph{wonderful compactification} of the complement of a hyperplane arrangement.
Surprisingly, for any matroid $M$, the Chow ring $A^*(\Sigma_M)$ has many of the properties of the cohomology ring of a smooth projective variety. 

We say a fan $\Sigma$ is \emph{Lefschetz} if its Chow ring $A^*(\Sigma)$  satisfies Poincar\'e duality, the hard Lefschetz theorem, and the Hodge-Riemann relations, and the star of any face in $\Sigma$ also has these properties. For a precise definition, see \cite{ArdilaDenhamHuh}.
Huh \cite{Huh}, Huh and Katz \cite{HuhKatz}, and Adiprasito, Huh, and Katz \cite{AdiprasitoHuhKatz} developed the first steps in the Hodge theory of matroids:

\begin{theorem} \label{thm:AHK} \cite{AdiprasitoHuhKatz}
The Bergman fan of a matroid is Lefschetz.
\end{theorem}

The inspiration for this theorem is geometric, coming from the Grothendieck standard conjectures on algebraic cycles. The statement and proof are combinatorial.

Instead of giving a complete definition of Lefschetz fans here, we focus on a comparatively small but powerful consequence. The Chow ring $A^*(\Sigma_M)$ is graded of degree $r-1$, and there is an isomorphism $\deg: A^{r-1} \rightarrow \R$ characterized by the property that $\deg(F_1\cdots F_{r-1})=1$ for any complete flag $F_1 \subsetneq \cdots \subsetneq F_{r-1}$ of proper flats.

Consider the \emph{ample cone} $K(\Sigma_M) \subset A^1(\Sigma_M)$ given by the piecewise linear functions on the Bergman fan $\Sigma_M$ that are strictly convex around every cone. In Brion's presentation, 
$
{K}(\Sigma_M) = \{ \sum_{F \textrm{ flat}} c_F x_F \, : \, 
c:2^E \rightarrow \R
 \, \textrm{ strictly submodular}\}
 $
The Hodge-Riemann relations imply that for any ample classes $L_1, \ldots, L_{r-3},a,b \in K(\Sigma_M)$, if we write $L=L_1 \cdots L_{r-3}$,  
\begin{equation}\label{eq:logconcave}
\deg(La^2)
\deg(Lb^2)
\leq
\deg(Lab)^2.
\end{equation}
By continuity, this property also holds for \emph{nef classes}, \emph{i.e.}, classes in the closure $\overline{K}(\Sigma_M)$.

Combining these ingredients, Adiprasito, Huh, and Katz \cite{AdiprasitoHuhKatz} considered the elements 
\[
\alpha \coloneq \alpha_i = \sum_{F \ni i} x_F, \qquad 
\beta \coloneq \beta_i = \sum_{F \not\owns i} x_F, 
\]
of the Chow ring $A^*(\Sigma_M)$, which are independent of $i$ and  lie in the nef cone $\overline{K}(\Sigma_M)$. An algebraic combinatorial computation in $A^*(\Sigma_M)$ shows that
\begin{equation}\label{eq:degalphabeta}
\deg(\alpha^k \beta^{r-1-k}) =
 (-1)^{r-1-k}\left( \textrm{coefficient of } q^k \textrm{ in } \frac{\chi_M(q)}{q-1}\right).
\end{equation}
As $k$ varies, this sequence of degrees is log-concave by  \eqref{eq:logconcave}. In turn, by elementary arguments \cite{Brylawski, JuhnkeKubitzkeLe, Lenz}, this  implies the following theorems, which were conjectured by Rota, Heron, Mason, and Welsh in the 1970s and 1980s.

\begin{theorem}  \label{thm:AHK2} \cite{AdiprasitoHuhKatz} For any matroid $M$ the following sequences, defined in Section \ref{sec:invariants}, are unimodal, log-concave, and flawless:

\noindent $\bullet$ 
the $f$-vector $f(\II(M))$ of the independence complex of $M$, and 

\noindent $\bullet$ 
the $f$-vector $f(\BC_<(M))$ of the broken circuit complex of $M$.
\end{theorem}

The latter is the sequence of absolute values of the coefficients of $\chi_M(q)$.

\subsection{Chern-Schwartz-MacPherson cycles of matroids} \label{sec:CSM}  

\emph{Chern-Schwartz-MacPherson} (CSM) cycles generalize the Chern class of a tangent bundle to the setting of possibly singular or non-compact complex algebraic varieties. When $\AA$ is a complex hyperplane arrangement, 
L\'opez de Medrano, Rinc\'on, and Shaw \cite{LRS} computed the CSM class of the complement $\C^E \setminus \AA$ of $\AA$  in its wonderful compactification $W_\AA$  in terms of the matroid $M = M(\AA)$.

A $k$-dimensional \emph{Minkowski weight} on a fan $\Sigma$ is a choice of weights $w(\sigma)$ for each $k$-dimensional face $\sigma$ of $\Sigma$ satisfying the balancing condition \eqref{eq:balancing} for every $(k-1)$-dimensional face $\tau$ of $\Sigma$. We write $\MW_k(\Sigma)$ for the additive group of $k$-dimensional Minkowski weights and $\MW(\Sigma) = \bigoplus_{k \geq 0} \MW_k(\Sigma)$. This is dual to the Chow ring of $\Sigma$ in the following sense. Fulton and Sturmfels \cite{FultonSturmfels} showed that $\MW_k(\Sigma) \cong \Hom(A^k(\Sigma),\R).$
The product in $A(\Sigma)$ then gives $\MW(\Sigma)$ the structure of an $A(\Sigma)-$module, and $\MW(\Sigma) \cong \Hom(A(\Sigma),\R)$ as modules. For details, see for example \cite[Section 3.1]{ArdilaDenhamHuh}.

The \emph{$k$-th CSM cycle of a matroid $M$} 
is the $k$-skeleton of the Bergman fan $\Sigma_M$ with weights 
\[
w(\sigma_\FF) \coloneq (-1)^{r-k} \prod_{i=0}^k \beta(M[F_i, F_{i+1}]) \qquad \textrm{ for } \FF = \{\emptyset = F_0 \subset F_1 \subset \cdots \subset F_k \subset F_{k+1} \subset E\}, 
\]
where $\beta(M[F_i, F_{i+1}])$ is the beta invariant of the minor $M[F_i, F_{i+1}]$ \cite{LRS}. 
For any matroid, the fan above satisfies the balancing condition \eqref{eq:balancing}, giving a Minkowski weight.

When $M$ is the matroid of a complex hyperplane arrangement $\AA$,  the (geometric) CSM class of the wonderful compactification $W_{\AA}$ is given by the (combinatorial) CSM cycles of $M$. 
The above construction makes sense for arbitrary matroids, and further, it defines the \emph{CSM cycles of tropical manifolds}. 

As shown in \cite{Aluffi, LRS}, the tropical degrees of the CSM cycles of a matroid $M$ are the entries of the $h$-vector of the broken circuit complex of $M$:
\begin{equation} \label{eq:degcsm}
\deg(\csm_k(M)) = \text{coefficient of $q^k$ in } \chi_M(q+1).
\end{equation}

\section{Geometric model 3: conormal fans}
\label{sec:geom3}

Motivated by Lagrangian geometry, Ardila, Denham, and Huh \cite{ArdilaDenhamHuh}  introduced a third polyhedral model that enriches the geometry of matroids and 
 leads to stronger inequalities for matroid invariants. In this section we assume that the matroid $M$ has no loops or coloops.
 A \emph{biflag} $(\FF,\GG)$ of $M$ consists of flags
 $\FF = \{F_1 \subseteq \cdots \subseteq F_l\}$ and $\GG = \{G_1 \supseteq \cdots \supseteq G_l\}$
 of nonempty flats of $M$ and $M^\perp$, respectively, such that
\[
\bigcap_{i=1}^l(F_i \cup G_i) = E, \qquad
\bigcup_{i=1}^l(F_i \cap G_i) \neq E.
\]
All maximal biflags have length $n-2$.

\begin{definition} \cite{ArdilaDenhamHuh}
The \emph{conormal fan} $\Sigma_{M,M^\perp}$ of a matroid $M$ is the polyhedral fan in $\R^E/\left< \e_E \right> \times \R^E/\left< \f_E \right>$ consisting of the cones
\[
\sigma_{\FF,\GG} \coloneq \cone \{\e_{F_i}+\f_{G_i} \, : \, 1 \leq i \leq l\} 
\qquad \text{for each biflag $(\FF,\GG)$}. 
\]
Here $\{\e_i : i \in E\}$ and $\{\f_i : i \in E\}$ are  the standard bases for two copies of $\R^E$. 
\end{definition}

This is a simplicial fan whose support is the product 
$\Trop M \times \Trop M^\perp$.

\subsection{The conormal Chow ring and combinatorial Hodge theory}

A \emph{biflat} $(F,G)$ of $M$ is a biflag of length $1$. It 
consists of flats $F, G \neq \emptyset$ of $M, M^\perp$ respectively, not both equal to $E$, such that $F \cup G = E$.
Consider the polynomial ring with a variable $x_{F,G}$ for each biflat $(F,G)$. 
For a set $(\FF,\GG)$ of distinct biflats $(F_1,G_1), \ldots, (F_l,G_l)$ write $x_{\FF,\GG} = x_{F_1, G_1} \cdots x_{F_l,G_l}$. For $i \in E$ let
\[
\gamma_i \coloneq \sum_{F \ni i \atop F \neq E} x_{F,G}, \quad
\gamma'_i \coloneq \sum_{G \ni i \atop G \neq E} x_{F,G}, \quad
\delta_i \coloneq \sum_{F \cap G \ni i} x_{F,G}.
\]
The \emph{Chow ring of the conormal fan} of $M$ is
\[
A^*(\Sigma_{M,M^\perp}) \coloneq \R[x_{F,G}]/(I_{M,M^\perp} + J_{M,M^\perp}),
\]
where
$
I_{M,M^\perp} = \left< x_{\FF,\GG} \, : \, (\FF ,\GG) \textrm{ is not a biflag }\right>$, $
J_{M, M^\perp} = \left<\gamma_i-\gamma_j, \gamma'_i-\gamma'_j \, : \, i,j \in E \right>.
$
The elements $\gamma \coloneq \gamma_i, \, \gamma' \coloneq \gamma'_i,$ and $\delta \coloneq \delta_i$ of $A^1(\Sigma_{M,M^\perp})$ are independent of $i$.

Ardila, Denham, and Huh \cite{ArdilaDenhamHuh} showed the conormal analog of Theorem \ref{thm:AHK}.

\begin{theorem} \label{thm:ArdilaDenhamHuh} \cite{ArdilaDenhamHuh}
The conormal fan of a matroid is Lefschetz.
\end{theorem}

Since $|\Sigma_{M,M^\perp}| = \Trop M \times \Trop M^\perp = |\Sigma_M| \times |\Sigma_{M^\perp}|$
and the product of Lefschetz fans is Lefschetz, the key step in the proof of Theorem \ref{thm:ArdilaDenhamHuh} is the general result
that the Lefschetz property of a simplicial fan $\Sigma$ depends only on the support $|\Sigma|$:

\begin{theorem} \label{thm:ArdilaDenhamHuh2} \cite{ArdilaDenhamHuh}
If two simplicial fans $\Sigma_1$ and $\Sigma_2$ have the same support $|\Sigma_1| = |\Sigma_2|$, then $\Sigma_1$ is Lefschetz if and only if $\Sigma_2$ is Lefschetz.
\end{theorem}

This Chow ring $A^*(\Sigma_{M,M^\perp})$ has degree $n-2$, and there is a unique isomorphism $\deg: A^{n-2} \rightarrow \R$ characterized by the property that $\deg(x_{\FF,\GG})=1$ for any maximal biflag $(\FF, \GG)$ of $M$. The log-concavity inequality \eqref{eq:logconcave} holds in the ample cone ${K}(\Sigma_{M,M^\perp})$ of the conormal fan, and hence in the nef cone $\overline{K}(\Sigma_{M,M^\perp})$ as well.

\subsection{Lagrangian interpretation of CSM classes}

We return to the Chern-Schwartz-MacPherson classes of Section \ref{sec:CSM}.
Schwartz's and MacPherson's constructions \cite{Schwartz, MacPherson} of $\csm$ for a complex algebraic variety $X$ are rather subtle. Sabbah  \cite{Sabbah} later observed that CSM 
classes can be interpreted more simply as ``shadows'' of the characteristic cycles in the cotangent bundle $T^*X$.

Similarly, the CSM cycles of a matroid $M$ are combinatorially intricate fans supported on the Bergman fan $\Sigma_M$.
We prove that they are ``shadows'' of  much simpler cycles of the conormal fan $\Sigma_{M, M^\perp}$. There is a natural projection map $\pi: \Sigma_{M,M^\perp} \rightarrow \Sigma_M$ which gives a pushforward map $\pi_{*}\colon \MW_k(\Sigma_{M, M^\perp}) \longrightarrow \MW_k(\Sigma_{M})$. We have:
\begin{theorem} \cite{ArdilaDenhamHuh} \label{thm:CSMTheorem}
If $M$ has no loops and no coloops, we have 
 \[
  \csm_k(M)=(-1)^{r-k} \pi_{*} (\delta^{n-k-1} \cap 1_{M,M^\perp} ) \qquad  \text{for $0 \le k \le r$},
 \]
where  $1_{M,M^\perp}$ is the top-dimensional constant Minkowski weight $1$ on the conormal fan. 
 \end{theorem}

\subsection{Unimodality, log-concavity, and flawlessness}

Applying the projection formula to Theorem \ref{thm:CSMTheorem} and \eqref{eq:degcsm}, we then express the $h$-vector of the broken circuit complex of $M$
 in the intersection theory of the conormal fan:
\begin{equation} \label{eq:deggammadelta}
\deg( \gamma^k\, \delta^{n-k-1}) = (-1)^{r-k}(\text{coefficient of $q^k$ in } \chi_M(q+1)).
\end{equation}
We give an alternative proof of this identity that does not rely on CSM classes in \cite{ArdilaDenhamHuh2}, through a careful study of the \emph{Lagrangian combinatorics of matroids}. 

The classes $\gamma$ and $\delta$ are nef, so the log-concavity inequalities \eqref{eq:logconcave} apply to the sequence of degrees in \eqref{eq:deggammadelta}.
This implies the following strengthening of Theorem \ref{thm:AHK2}, parts of which were originally conjectured by Brylawski, Dawson, and Colbourn in the early 80s \cite{Brylawski, Colbourn, Dawson} and left open in Huh's 2018 ICM paper \cite{HuhICM}.

\begin{theorem}  \label{thm:ArdilaDenhamHuh3} \cite{ArdilaDenhamHuh} 
For any matroid $M$ the following sequences, defined in Section \ref{sec:invariants}, are unimodal, log-concave, and flawless:

\noindent $\bullet$ 
the $h$-vector $h(\II(M))$ of the independence complex of $M$, and 

\noindent $\bullet$ 
the $h$-vector $h(\BC_<(M))$ of the broken circuit complex of $M$.
\end{theorem}

The most difficult part of Theorem \ref{thm:ArdilaDenhamHuh3} is the log-concavity of 
$h(\BC_<(M))$ \eqref{eq:deggammadelta}. 
The remaining parts follow from it by elementary arguments.

In 1977, 
Stanley \cite{Stanleyh-vector} conjectured that $h(\II(M))$ is the $f$-vector of a \emph{pure multicomplex}: a set of monomials such that if $m' \in X$ and $m|m'$ then $m \in X$, and the maximal monomials in $X$ have the same degree. This conjecture has been proved in rather different ways for various families of matroids, \emph{e.g.} \cite{Merinocographic, Ohcotransversal}, but remains wide open in general.

\subsection{Lagrangian combinatorics, bipermutahedra, and harmonic polytopes}

A subtle technical issue in the proofs of Theorems \ref{thm:AHK2} and \ref{thm:ArdilaDenhamHuh3} leads to some combinatorial constructions of independent interest. 
For a Lefschetz fan $\Sigma$, the log-concavity inequalities \eqref{eq:logconcave} hold inside the open cone $K(\Sigma)$, corresponding to piecewise linear functions on $\Sigma$ that are \textbf{strictly} convex around every cone.
In light of \eqref{eq:degalphabeta} and \eqref{eq:deggammadelta}, we wish to apply these inequalities to the classes $\alpha, \beta$ in $A^1(\Sigma_M)$ and $\gamma, \delta$ in $A^1(\Sigma_{M, M^\perp})$, which are \textbf{weakly} convex locally; that is, they lie in the closed cone $\overline{K}(\Sigma)$. 
Continuity will guarantee that the log-concavity inequalities in $K(\Sigma)$ will still hold in the closure $\overline{K}(\Sigma)$ -- \textbf{if} the open cone $K(\Sigma)$ is non-empty. This is not a trivial condition.

A complete polyhedral fan $\Sigma$ is \emph{projective} if it is the normal fan $\Sigma = \mathcal{N}(P)$ of a polytope $P$; each such polytope produces a strictly convex piecewise linear function on the fan, namely, $w \mapsto \max_{p \in P} w(p)$. Under this correspondence, we can think of $K(\Sigma)$ as the space of polytopes whose normal fan is $\Sigma$. 
A non-complete fan $\Sigma$ is \emph{quasi-projective} if it is a subfan of a projective fan $\mathcal{N}(P)$; this guarantees that $K(\Sigma) \neq \emptyset$ as well.

By construction, the Bergman fan $\Sigma_M$ of every matroid $M$ on $E$ is a subfan of the braid arrangement $\Sigma_E$; this is the normal fan of the permutahedron $\Pi_E$, which is simple. This guarantees $K(\Sigma_M) \neq \emptyset$, as required in the proof of Theorem \ref{thm:AHK2}. Similarly, we had to construct a simple polytope, called the \emph{bipermutahedron} $\Pi_{E,E}$, whose normal fan contains the conormal fan $\Sigma_{M,M^\perp}$ of any matroid $M$ on $E$.
This guarantees $K(\Sigma_{M.M^\perp}) \neq \emptyset$, as required in the proof of Theorem \ref{thm:ArdilaDenhamHuh3}. Let us briefly discuss this polytope.

\smallskip

\noindent \textbf{\textsf{The harmonic polytope.}} A \emph{bisubset} of $E$, denoted $S|T \sqsubset E$, consists of subsets $S,T \neq \emptyset$ of $E$, not both equal to $E$, with $S \cup T = E$. Ardila and Escobar \cite{ArdilaEscobar} studied the \emph{harmonic polytope} $H_{E,E}$, given by
\[
\sum_{e \in E} x_e = \sum_{e \in E} y_e = \frac{n(n+1)}2 + 1, \quad
\sum_{s \in S} x_s + \sum_{t \in T} y_t \geq \frac{|S|(|S|+1) + |T|(|T|+1)}2  + 1   \text{ for } S|T \sqsubset E. 
\]
It has $3^n-3$ facets and $(n!)^2(1+\frac12+\cdots+\frac1n)$ vertices. Its volume is a weighted sum of the degrees of the toric varieties associated to all 
connected bipartite graphs with $n$ edges.

The harmonic polytope is the minimal polytope whose normal fan contains all conormal fans as subfans: every such polytope contains $H_{E,E}$ as a Minkowski summand. Thus the harmonic polytope arises very naturally in this setting, but it has the disadvantage of not being simple, which makes it difficult to work with its Chow ring.

\smallskip

\noindent \textbf{\textsf{The bipermutahedron.}} 
The \emph{bipermutahedron} $\Pi_{E,E} \subset \R^E \times \R^E$ is given by: 
\[
\sum_{e \in E} x_e =  \sum_{e \in [n]} y_e = 0, \quad  \sum_{s \in S} x_s + \sum_{t \in T} y_t \geq - \, \big(|S|+|S-T|\big)  \big(|T|+|T-S|\big)  \text{ for } S|T \sqsubset E. 
\]
It was constructed by Ardila, Denham, and Huh \cite{ArdilaDenhamHuh} and further studied in \cite{Ardilabiperm}. 
This is the most elegant simple polytope that we know which has the harmonic polytope as a Minkowski summand, so its normal fan contains all conormal fans. Its normal fan is simplicial, so we can use Brion and Billera's descriptions of its Chow ring. Finally, it has an elegant combinatorial structure that allows us to prove Theorem \ref{thm:CSMTheorem} and  \eqref{eq:deggammadelta} in \cite{ArdilaDenhamHuh}.

The bipermutahedron also has  $3^n-3$ facets,
 and it has $(2n)!/2^n$ vertices corresponding to the \emph{bipermutations of} $E$. It is combinatorially isomorphic to a unimodular triangulation of the product of $n$ unit triangles, and Ehrhart theory then gives a simple formula for its $h$-vector.
In analogy with the permutahedron, this $h$-vector counts the bipermutations according to their number of descents, and it is also log-concave \cite{Ardilabiperm}.

\section{Geometry and geometries: the future}
\label{sec:future}

We have centered our discussion on three related geometric models of a combinatorial geometry and their consequences within and outside of matroid theory. There are certainly many other such models -- some already known, some yet to be discovered. What these constructions share is their deep connection with or analogy to natural geometric constructions associated to vector configurations or subspaces of a vector space.
This situation reminds us of a 50-year old prophetic remark of Bose, relayed by Kelly and Rota: 

\begin{quote}
\vspace{-.2cm}
{``We combinatorialists have much to gain from the study of algebraic geometry, if not by its applications to our field, at least by the analogies between the two subjects."} \cite{KellyRota}
\end{quote}
\vspace{-.2cm}

\noindent 
The recurring theme of developing discrete versions of geometric techniques is not born from a wish to avoid algebraic geometry; quite the opposite. Our goal is to develop the necessary tools to solve combinatorial and geometric problems when the current algebro-geometric technology is not sufficient.  The applications of this program are not only combinatorial. As Gelfand, Goresky, MacPherson, and Serganova wrote in 1987,

\vspace{-.2cm}
\begin{quote}
{``We believe that combinatorial methods will play an increasing role in the future of geometry and topology."}  \cite{GGMS}
\end{quote}
\vspace{-.2cm}

\noindent 
Today, these predictions ring true more than ever.

\begin{ack}
I warmly thank Gian-Carlo Rota, Richard Stanley, and Bernd Sturmfels for their very beautiful and very different lessons on combinatorial geometries, and my wonderful coauthors and students for the great fun we've had understanding matroids together. I am also grateful to
Carolina Benedetti,
Dusty Ross, 
Felipe Rinc\'on,
Graham Denham,
Johannes Rau,
June Huh,
Kris Shaw,
Laura Escobar,
Mariel Supina, 
Mario Sanchez, 
May-Li Khoe, and
Ra\'ul Penagui\~ao 
for their valuable feedback on a draft of this manuscript. 
\end{ack}

\begin{funding}
This work was partially supported by United States National Science Foundation grant  DMS-1855610 and Simons Fellowship 613384.
\end{funding}








\bibliographystyle{myamsplain}
\bibliography{references}

\def\cprime{$'$} \def\cprime{$'$}
\providecommand{\bysame}{\leavevmode\hbox to3em{\hrulefill}\thinspace}
\providecommand{\MR}{\relax\ifhmode\unskip\space\fi MR }
\providecommand{\MRhref}[2]{%
  \href{http://www.ams.org/mathscinet-getitem?mr=#1}{#2}
}
\providecommand{\href}[2]{#2}
\begin{thebibliography}{100}

\bibitem{AdiprasitoHuhKatz}
K.~Adiprasito, J.~Huh, and E.~Katz, \emph{Hodge theory for combinatorial
  geometries}, Ann. of Math. (2) \textbf{188} (2018), no.~2, 381--452.
  \MR{3862944}

\bibitem{AguiarArdila}
M.~Aguiar and F.~Ardila, \emph{Hopf monoids and generalized permutahedra},
  arXiv:1709.07504, Memoirs of the Amer. Math. Soc. (to appear).

\bibitem{Alessandrini}
D.~Alessandrini, \emph{Logarithmic limit sets of real semi-algebraic sets},
  Adv. Geom. \textbf{13} (2013), no.~1, 155--190. \MR{3011539}

\bibitem{AllermannRau}
L.~Allermann and J.~Rau, \emph{First steps in tropical intersection theory},
  Mathematische Zeitschrift \textbf{264} (2010), no.~3, 633--670.

\bibitem{Aluffi}
P.~Aluffi, \emph{Grothendieck classes and {C}hern classes of hyperplane
  arrangements}, Int. Math. Res. Not. IMRN (2013), no.~8, 1873--1900.
  \MR{3047491}

\bibitem{ArdilaThesis}
F.~Ardila, \emph{Enumerative and algebraic aspects of matroids and hyperplane
  arrangements}, Ph.D. thesis, Massachusetts Institute of Technology, 2003.

\bibitem{ArdilaTutte}
\bysame, \emph{Computing the {T}utte polynomial of a hyperplane arrangement},
  Pacific J. Math. \textbf{230} (2007), no.~1, 1--26. \MR{2318445
  (2008g:52034)}

\bibitem{Ardilasurvey}
\bysame, \emph{Algebraic and geometric methods in enumerative combinatorics},
  Handbook of enumerative combinatorics, Discrete Math. Appl. (Boca Raton), CRC
  Press, Boca Raton, FL, 2015, pp.~3--172. \MR{3409342}

\bibitem{Ardilabiperm}
\bysame, \emph{The bipermutahedron}, arXiv:2008.02295 (2020).

\bibitem{ArdilaTuttesurvey}
\bysame, \emph{Tutte polynomials of hyperplane arrangements and the finite
  field method}, Handbook of the Tutte polynomial and related topics, CRC
  Press, Boca Raton, FL, to appear.

\bibitem{ArdilaBeckMcWhirter}
F.~Ardila, M.~Beck, and J.~McWhirter, \emph{The arithmetic of {C}oxeter
  permutahedra}, Rev. Acad. Colombiana Cienc. Exact. F\'{\i}s. Natur.
  \textbf{44} (2020), no.~173, 1152--1166. \MR{4236163}

\bibitem{ArdilaBenedettiDoker}
F.~Ardila, C.~Benedetti, and J.~Doker, \emph{Matroid polytopes and their
  volumes}, Discrete Comput. Geom. \textbf{43} (2010), no.~4, 841--854.
  \MR{2610473 (2012b:52026)}

\bibitem{ArdilaCastilloEurPostnikov}
F.~Ardila, F.~Castillo, C.~Eur, and A.~Postnikov, \emph{Coxeter submodular
  functions and deformations of {C}oxeter permutahedra}, Adv. Math.
  \textbf{365} (2020), 107039, 36. \MR{4064768}

\bibitem{ArdilaCastilloHenley}
F.~Ardila, F.~Castillo, and M.~Henley, \emph{The arithmetic {T}utte polynomials
  of the classical root systems}, International Mathematics Research Notices
  \textbf{2015} (2014), no.~12, 3830--3877.

\bibitem{ArdilaDenhamHuh}
F.~Ardila, G.~Denham, and J.~Huh, \emph{Lagrangian geometry of matroids},
  arXiv:2004.13116 (2020).

\bibitem{ArdilaDenhamHuh2}
\bysame, \emph{Lagrangian combinatorics of matroids}, arXiv:2109.11565 (2021).

\bibitem{ArdilaEscobar}
F.~Ardila and L.~Escobar, \emph{The harmonic polytope}, Selecta Math. (N.S.)
  \textbf{27} (2021), no.~5, Paper No. 91, 31. \MR{4314249}

\bibitem{ArdilaFinkRincon}
F.~Ardila, A.~Fink, and F.~Rinc{\'o}n, \emph{Valuations for matroid polytope
  subdivisions}, Canad. J. Math. \textbf{62} (2010), no.~6, 1228--1245.
  \MR{2760656 (2012c:05075)}

\bibitem{ArdilaKlivansWilliams}
F.~Ardila, C.~Klivans, and L.~Williams, \emph{The positive {B}ergman complex of
  an oriented matroid}, European J. Combin. \textbf{27} (2006), no.~4,
  577--591. \MR{2215218}

\bibitem{ArdilaKlivans}
F.~Ardila and C.~J. Klivans, \emph{The {B}ergman complex of a matroid and
  phylogenetic trees}, J. Combin. Theory Ser. B \textbf{96} (2006), no.~1,
  38--49. \MR{2185977 (2006i:05034)}

\bibitem{ArdilaReinerWilliams}
F.~Ardila, V.~Reiner, and L.~Williams, \emph{Bergman complexes, {C}oxeter
  arrangements, and graph associahedra}, S\'{e}m. Lothar. Combin. \textbf{54A}
  (2005/07), Art. B54Aj, 25. \MR{2317682}

\bibitem{ArdilaRinconWilliams2}
F.~Ardila, F.~Rinc{\'o}n, and L.~K. Williams, \emph{Positively oriented
  matroids are realizable}, Journal of the European Mathematical Society
  \textbf{19} (2017), no.~3, 815--833.

\bibitem{ArdilaSanchez}
F.~Ardila and M.~Sanchez, \emph{Valuations and the {H}opf monoid of generalized
  permutahedra}, arXiv:2010.11178 (2020).

\bibitem{Athanasiadis}
C.~A. Athanasiadis, \emph{Characteristic polynomials of subspace arrangements
  and finite fields}, Adv. Math. \textbf{122} (1996), no.~2, 193--233.
  \MR{1409420 (97k:52012)}

\bibitem{BallaOlarte}
G.~Balla and J.~A. Olarte, \emph{The tropical symplectic {G}rassmannian},
  International Mathematics Research Notices (2021), rnab267.

\bibitem{Bergman}
G.~M. Bergman, \emph{The logarithmic limit-set of an algebraic variety}, Trans.
  Amer. Math. Soc. \textbf{157} (1971), 459--469. \MR{280489}

\bibitem{BieriGroves}
R.~Bieri and J.~R.~J. Groves, \emph{The geometry of the set of characters
  induced by valuations}, J. Reine Angew. Math. \textbf{347} (1984), 168--195.
  \MR{733052}

\bibitem{Billera}
L.~J. Billera, \emph{The algebra of continuous piecewise polynomials}, Adv.
  Math. \textbf{76} (1989), no.~2, 170--183. \MR{1013666}

\bibitem{Bjornermatroids}
A.~Bj\"{o}rner, \emph{The homology and shellability of matroids and geometric
  lattices}, Matroid applications, Encyclopedia Math. Appl., vol.~40, Cambridge
  Univ. Press, Cambridge, 1992, pp.~226--283. \MR{1165544}

\bibitem{BjornerEkedahl}
A.~Bj{\"o}rner and T.~Ekedahl, \emph{Subspace arrangements over finite fields:
  cohomological and enumerative aspects}, Advances in Mathematics \textbf{129}
  (1997), no.~2, 159--187.

\bibitem{OMs}
A.~Bj{\"o}rner, M.~Las~Vergnas, B.~Sturmfels, N.~White, and G.~M. Ziegler,
  \emph{Oriented matroids}, second ed., Encyclopedia of Mathematics and its
  Applications, vol.~46, Cambridge University Press, Cambridge, 1999.
  \MR{1744046 (2000j:52016)}

\bibitem{Coxetermatroids}
A.~V. Borovik, I.~M. Gelfand, and N.~White, \emph{Coxeter matroids}, Progress
  in Mathematics, vol. 216, Birkh\"auser Boston, Inc., Boston, MA, 2003.
  \MR{1989953 (2004i:05028)}

\bibitem{BrandenMoci}
P.~Br{\"a}nd{\'e}n and L.~Moci, \emph{The multivariate arithmetic {T}utte
  polynomial}, Transactions of the American Mathematical Society \textbf{366}
  (2014), no.~10, 5523--5540.

\bibitem{BrandtEurZhang}
M.~Brandt, C.~Eur, and L.~Zhang, \emph{Tropical flag varieties}, Advances in
  Mathematics \textbf{384} (2021), 107695.

\bibitem{Brion}
M.~Brion, \emph{Piecewise polynomial functions, convex polytopes and
  enumerative geometry}, Parameter spaces ({W}arsaw, 1994), Banach Center
  Publ., vol.~36, Polish Acad. Sci. Inst. Math., Warsaw, 1996, pp.~25--44.
  \MR{1481477}

\bibitem{Brylawski}
T.~Brylawski, \emph{A combinatorial model for series-parallel networks}, Trans.
  Amer. Math. Soc. \textbf{154} (1971), 1--22.

\bibitem{Colbourn}
C.~J. Colbourn, \emph{The combinatorics of network reliability}, International
  Series of Monographs on Computer Science, The Clarendon Press, Oxford
  University Press, New York, 1987. \MR{902584}

\bibitem{CrapoRota}
H.~H. Crapo and G.-C. Rota, \emph{On the foundations of combinatorial theory:
  {C}ombinatorial geometries}, preliminary ed., The M.I.T. Press, Cambridge,
  Mass.-London, 1970. \MR{0290980 (45 \#74)}

\bibitem{daSilva}
I.~P.~F. Da~Silva, \emph{Quelques propri{\'e}t{\'e}s des matroides
  orient{\'e}s}, Ph.D. thesis, Universit{\'e} Pierre-et-Marie-Curie [- Paris
  VI, 1987.

\bibitem{D'AdderioMoci.Ehrhart}
M.~D'Adderio and L.~Moci, \emph{{E}hrhart polynomial and arithmetic {T}utte
  polynomial}, European J. Combin. \textbf{33} (2012), no.~7, 1479--1483.
  \MR{2923464}

\bibitem{DahmenMicchelli}
W.~Dahmen and C.~A. Micchelli, \emph{On the local linear independence of
  translates of a box spline}, Studia Math. \textbf{82} (1985), no.~3,
  243--263. \MR{825481 (87k:41008)}

\bibitem{Dawson}
J.~E. Dawson, \emph{A collection of sets related to the {T}utte polynomial of a
  matroid}, Graph theory, {S}ingapore 1983, Lecture Notes in Math., vol. 1073,
  Springer, Berlin, 1984, pp.~193--204. \MR{761018}

\bibitem{DeConciniProcesi.toric}
C.~De~Concini and C.~Procesi, \emph{On the geometry of toric arrangements},
  Transform. Groups \textbf{10} (2005), no.~3-4, 387--422. \MR{2183118
  (2006m:32027)}

\bibitem{DeConciniProcesi.Tutte}
\bysame, \emph{The zonotope of a root system}, Transform. Groups \textbf{13}
  (2008), no.~3-4, 507--526. \MR{2452603}

\bibitem{DeConciniProcesiVergne}
C.~De~Concini, C.~Procesi, and M.~Vergne, \emph{Vector partition functions and
  index of transversally elliptic operators}, Transform. Groups \textbf{15}
  (2010), no.~4, 775--811. \MR{2753257}

\bibitem{DerksenFink}
H.~Derksen and A.~Fink, \emph{Valuative invariants for polymatroids}, Adv.
  Math. \textbf{225} (2010), no.~4, 1840--1892. \MR{2680193 (2011k:52016)}

\bibitem{DressWenzel}
A.~W.~M. Dress and W.~Wenzel, \emph{Valuated matroids}, Adv. Math. \textbf{93}
  (1992), no.~2, 214--250. \MR{1164708 (93h:05045)}

\bibitem{Edmonds}
J.~Edmonds, \emph{Submodular functions, matroids, and certain polyhedra},
  Combinatorial structures and their applications (1970), 69--87.

\bibitem{EhrenborgReaddySlone}
R.~Ehrenborg, M.~Readdy, and M.~Slone, \emph{Affine and toric hyperplane
  arrangements}, Discrete Comput. Geom. \textbf{41} (2009), no.~4, 481--512.
  \MR{2496314 (2010a:52035)}

\bibitem{Ehrhart}
E.~Ehrhart, \emph{Sur les poly\`edres rationnels homoth\'{e}tiques \`a {$n$}
  dimensions}, C. R. Acad. Sci. Paris \textbf{254} (1962), 616--618.
  \MR{130860}

\bibitem{EPW}
B.~Elias, N.~Proudfoot, and M.~Wakefield, \emph{The {K}azhdan-{L}usztig
  polynomial of a matroid}, Adv. Math. \textbf{299} (2016), 36--70.
  \MR{3519463}

\bibitem{Eurvolume}
C.~Eur, \emph{Divisors on matroids and their volumes}, J. Combin. Theory Ser. A
  \textbf{169} (2020), 105135, 31. \MR{4011081}

\bibitem{EurSanchezSupina}
C.~Eur, M.~Sanchez, and M.~Supina, \emph{The universal valuation of {C}oxeter
  matroids}, Bull. Lond. Math. Soc. \textbf{53} (2021), no.~3, 798--819.
  \MR{4275090}

\bibitem{FeichtnerYuzvinsky}
E.~M. Feichtner and S.~Yuzvinsky, \emph{Chow rings of toric varieties defined
  by atomic lattices}, Inventiones Mathematicae \textbf{155} (2004), no.~3,
  515--536.

\bibitem{Fink}
A.~Fink, \emph{Tropical cycles and chow polytopes}, Beitr{\"a}ge zur Algebra
  und Geometrie/Contributions to Algebra and Geometry \textbf{54} (2013),
  no.~1, 13--40.

\bibitem{FinkSpeyer}
A.~Fink and D.~E. Speyer, \emph{K-classes for matroids and equivariant
  localization}, Duke Mathematical Journal \textbf{161} (2012), no.~14,
  2699--2723.

\bibitem{FultonSturmfels}
W.~Fulton and B.~Sturmfels, \emph{Intersection theory on toric varieties},
  Topology \textbf{36} (1997), no.~2, 335--353. \MR{1415592 (97h:14070)}

\bibitem{Geldon}
T.~W. Geldon, \emph{Computing the {T}utte polynomial of hyperplane
  arrangements}, ProQuest LLC, Ann Arbor, MI, 2009, Thesis (Ph.D.)--The
  University of Texas at Austin. \MR{2713747}

\bibitem{GGMS}
I.~M. Gel{\cprime}fand, R.~M. Goresky, R.~D. MacPherson, and V.~V. Serganova,
  \emph{Combinatorial geometries, convex polyhedra, and {S}chubert cells}, Adv.
  in Math. \textbf{63} (1987), no.~3, 301--316. \MR{877789 (88f:14045)}

\bibitem{Greeneweight}
C.~Greene, \emph{Weight enumeration and the geometry of linear codes}, Studies
  in Applied Mathematics \textbf{55} (1976), no.~2, 119--128.

\bibitem{HackingKeelTevelev}
P.~Hacking, S.~Keel, and J.~Tevelev, \emph{Compactification of the moduli space
  of hyperplane arrangements}, J. Algebraic Geom. \textbf{15} (2006), no.~4,
  657--680. \MR{2237265 (2007j:14016)}

\bibitem{Huh}
J.~Huh, \emph{Milnor numbers of projective hypersurfaces and the chromatic
  polynomial of graphs}, J. Amer. Math. Soc. \textbf{25} (2012), no.~3,
  907--927. \MR{2904577}

\bibitem{HuhICM}
\bysame, \emph{Combinatorial applications of the {H}odge-{R}iemann relations},
  Proceedings of the {I}nternational {C}ongress of {M}athematicians---{R}io de
  {J}aneiro 2018. {V}ol. {IV}. {I}nvited lectures, World Sci. Publ.,
  Hackensack, NJ, 2018, pp.~3093--3111. \MR{3966524}

\bibitem{HuhKatz}
J.~Huh and E.~Katz, \emph{Log-concavity of characteristic polynomials and the
  {B}ergman fan of matroids}, Math. Ann. \textbf{354} (2012), no.~3,
  1103--1116. \MR{2983081}

\bibitem{JKU}
D.~Jensen, M.~Kutler, and J.~Usatine, \emph{The motivic zeta functions of a
  matroid}, J. Lond. Math. Soc. (2) \textbf{103} (2021), no.~2, 604--632.
  \MR{4230913}

\bibitem{JoniRota}
S.~A. Joni and G.-C. Rota, \emph{Coalgebras and bialgebras in combinatorics},
  Umbral Calculus and Hopf Algebras (Norman, OK, 1978), Contemp. Math., vol.~6,
  Amer. Math. Soc., Providence, R.I., 1982, pp.~1--47.

\bibitem{JuhnkeKubitzkeLe}
M.~Juhnke-Kubitzke and D.~V. Le, \emph{Flawlessness of {$h$}-vectors of broken
  circuit complexes}, Int. Math. Res. Not. IMRN (2018), no.~5, 1347--1367.
  \MR{3801465}

\bibitem{Kapranov}
M.~M. Kapranov, \emph{Chow quotients of {G}rassmannians. {I}}, I. {M}.
  {G}el\cprime fand {S}eminar, Adv. Soviet Math., vol.~16, Amer. Math. Soc.,
  Providence, RI, 1993, pp.~29--110. \MR{1237834 (95g:14053)}

\bibitem{KellyRota}
D.~Kelly and G.-C. Rota, \emph{Some problems in combinatorial geometry}, A
  survey of combinatorial theory, Elsevier, 1973, pp.~309--312.

\bibitem{KnutsonTao}
A.~Knutson and T.~Tao, \emph{Puzzles and (equivariant) cohomology of
  {G}rassmannians}, Duke Mathematical Journal \textbf{119} (2003), no.~2,
  221--260.

\bibitem{Lafforgue}
L.~Lafforgue, \emph{Chirurgie des {G}rassmanniennes}, CRM Monograph Series,
  vol.~19, American Mathematical Society, Providence, RI, 2003. \MR{1976905
  (2004k:14085)}

\bibitem{Lason}
M.~Laso\'{n}, \emph{On the toric ideals of matroids of a fixed rank}, Selecta
  Math. (N.S.) \textbf{27} (2021), no.~2, Paper No. 18, 17. \MR{4236584}

\bibitem{Lenz}
M.~Lenz, \emph{The {$f$}-vector of a representable-matroid complex is
  log-concave}, Adv. in Appl. Math. \textbf{51} (2013), no.~5, 543--545.
  \MR{3118543}

\bibitem{LRS}
L.~L\'{o}pez~de Medrano, F.~Rinc\'{o}n, and K.~Shaw,
  \emph{Chern-{S}chwartz-{M}ac{P}herson cycles of matroids}, Proc. Lond. Math.
  Soc. (3) \textbf{120} (2020), no.~1, 1--27. \MR{3999674}

\bibitem{MaclaganSturmfels}
D.~Maclagan and B.~Sturmfels, \emph{Introduction to {T}ropical {G}eometry},
  Graduate Studies in Mathematics, vol. 161, American Mathematical Society,
  Providence, RI, 2015.

\bibitem{MacPherson}
R.~D. MacPherson, \emph{Chern classes for singular algebraic varieties}, Ann.
  of Math. (2) \textbf{100} (1974), 423--432. \MR{361141}

\bibitem{McMullen}
P.~McMullen, \emph{The polytope algebra}, Advances in Mathematics \textbf{78}
  (1989), no.~1, 76--130.

\bibitem{Merinocographic}
C.~Merino, \emph{The chip firing game and matroid complexes}, Discrete models:
  combinatorics, computation, and geometry ({P}aris, 2001), Discrete Math.
  Theor. Comput. Sci. Proc., AA, Maison Inform. Math. Discr\`et. (MIMD), Paris,
  2001, pp.~245--255 (electronic). \MR{1888777 (2002k:91048)}

\bibitem{Mikhalkin}
G.~Mikhalkin, \emph{Enumerative tropical algebraic geometry in {$\mathbb
  R^2$}}, J. Amer. Math. Soc. \textbf{18} (2005), no.~2, 313--377. \MR{2137980}

\bibitem{MikhalkinRau}
G.~Mikhalkin and J.~Rau, \emph{Tropical geometry}, 2018.

\bibitem{Moci.toric}
L.~Moci, \emph{A {T}utte polynomial for toric arrangements}, Trans. Amer. Math.
  Soc. \textbf{364} (2012), no.~2, 1067--1088. \MR{2846363}

\bibitem{Murota}
K.~Murota, \emph{Convexity and {S}teinitz's exchange property}, Adv. Math.
  \textbf{124} (1996), no.~2, 272--311. \MR{1424312 (97m:90070)}

\bibitem{Nakasawa}
T.~Nakasawa, \emph{Zur axiomatik der linearen abh{\"a}ngigkeit. i}, Science
  Reports of the Tokyo Bunrika Daigaku, Section A \textbf{2} (1935), no.~43,
  235--255.

\bibitem{Nelson}
P.~Nelson, \emph{Almost all matroids are nonrepresentable}, Bull. Lond. Math.
  Soc. \textbf{50} (2018), no.~2, 245--248. \MR{3830117}

\bibitem{Ohcotransversal}
S.~Oh, \emph{Generalized permutohedra, {$h$}-vectors of cotransversal matroids
  and pure {O}-sequences}, Electron. J. Combin. \textbf{20} (2013), no.~3,
  Paper 14, 14. \MR{3104512}

\bibitem{OrlikSolomon}
P.~Orlik and L.~Solomon, \emph{Combinatorics and topology of complements of
  hyperplanes}, Invent. Math. \textbf{56} (1980), no.~2, 167--189. \MR{558866
  (81e:32015)}

\bibitem{Oxley}
J.~Oxley, \emph{Matroid theory}, second ed., Oxford Graduate Texts in
  Mathematics, vol.~21, Oxford University Press, Oxford, 2011. \MR{2849819
  (2012k:05002)}

\bibitem{Postnikovgenperm}
A.~Postnikov, \emph{Permutohedra, associahedra, and beyond}, Int. Math. Res.
  Not. IMRN (2009), no.~6, 1026--1106. \MR{2487491 (2010g:05399)}

\bibitem{Randriamaro}
H.~Randriamaro, \emph{The {T}utte polynomial of symmetric hyperplane
  arrangements}, J. Knot Theory Ramifications \textbf{29} (2020), no.~3,
  2050004, 19. \MR{4101598}

\bibitem{RinconD}
F.~Rinc{\'o}n, \emph{Isotropical linear spaces and valuated delta-matroids},
  Journal of Combinatorial Theory, Series A \textbf{119} (2012), no.~1, 14--32.

\bibitem{RotaMobius}
G.-C. Rota, \emph{On the foundations of combinatorial theory. {I}. {T}heory of
  {M}\"{o}bius functions}, Z. Wahrscheinlichkeitstheorie und Verw. Gebiete
  \textbf{2} (1964), 340--368 (1964). \MR{174487}

\bibitem{RotaICM}
\bysame, \emph{Combinatorial theory, old and new}, Actes du {C}ongr\`es
  {I}nternational des {M}ath\'{e}maticiens ({N}ice, 1970), {T}ome 3,
  Gauthier-Villars / Paris, 1971, pp.~229--233. \MR{0505646}

\bibitem{Sabbah}
C.~Sabbah, \emph{Quelques remarques sur la g\'{e}om\'{e}trie des espaces
  conormaux}, Ast\'{e}risque (1985), no.~130, 161--192, Differential systems
  and singularities (Luminy, 1983). \MR{804052}

\bibitem{Saganfactors}
B.~E. Sagan, \emph{Why the characteristic polynomial factors}, Bull. Amer.
  Math. Soc. (N.S.) \textbf{36} (1999), no.~2, 113--133. \MR{1659875
  (2000a:06021)}

\bibitem{Schmitt}
W.~R. Schmitt, \emph{Antipodes and incidence coalgebras}, Journal of
  Combinatorial Theory, Series A \textbf{46} (1987), no.~2, 264--290.

\bibitem{Schwartz}
M.-H. Schwartz, \emph{Classes caract\'{e}ristiques d\'{e}finies par une
  stratification d'une vari\'{e}t\'{e} analytique complexe. {I}}, C. R. Acad.
  Sci. Paris \textbf{260} (1965), 3262--3264. \MR{212842}

\bibitem{SpeyerWilliams}
D.~Speyer and L.~Williams, \emph{The tropical totally positive {G}rassmannian},
  J. Algebraic Combin. \textbf{22} (2005), no.~2, 189--210. \MR{2164397}

\bibitem{Speyer1}
D.~E. Speyer, \emph{Tropical linear spaces}, SIAM J. Discrete Math. \textbf{22}
  (2008), no.~4, 1527--1558. \MR{2448909 (2009m:52023)}

\bibitem{Speyer2}
\bysame, \emph{A matroid invariant via the {$K$}-theory of the {G}rassmannian},
  Adv. Math. \textbf{221} (2009), no.~3, 882--913. \MR{2511042 (2010i:14114)}

\bibitem{Stanleyh-vector}
R.~P. Stanley, \emph{Cohen-{M}acaulay complexes}, Higher combinatorics ({P}roc.
  {NATO} {A}dvanced {S}tudy {I}nst., {B}erlin, 1976), Reidel, Dordrecht, 1977,
  pp.~51--62. NATO Adv. Study Inst. Ser., Ser. C: Math. and Phys. Sci., 31.
  \MR{0572989 (58 \#28010)}

\bibitem{StanleyLie}
\bysame, \emph{Unimodal sequences arising from {L}ie algebras}, Combinatorics,
  representation theory and statistical methods in groups, Lecture Notes in
  Pure and Appl. Math., vol.~57, Dekker, New York, 1980, pp.~127--136.
  \MR{588199}

\bibitem{StanleyICM}
\bysame, \emph{Combinatorial applications of the hard {L}efschetz theorem},
  Proceedings of the {I}nternational {C}ongress of {M}athematicians, {V}ol. 1,
  2 ({W}arsaw, 1983), PWN, Warsaw, 1984, pp.~447--453. \MR{804700}

\bibitem{Stanleyzonotope}
\bysame, \emph{A zonotope associated with graphical degree sequences}, Applied
  geometry and discrete mathematics, DIMACS Ser. Discrete Math. Theoret.
  Comput. Sci., vol.~4, Amer. Math. Soc., Providence, RI, 1991, pp.~555--570.
  \MR{1116376 (92k:52020)}

\bibitem{Sturmfelstoriceqs}
B.~Sturmfels, \emph{Equations defining toric varieties}, Algebraic
  geometry---{S}anta {C}ruz 1995, Proc. Sympos. Pure Math., vol.~62, Amer.
  Math. Soc., Providence, RI, 1997, pp.~437--449. \MR{1492542}

\bibitem{Sturmfelspolyeqs}
\bysame, \emph{Solving systems of polynomial equations}, CBMS Regional
  Conference Series in Mathematics, vol.~97, Published for the Conference Board
  of the Mathematical Sciences, Washington, DC; by the American Mathematical
  Society, Providence, RI, 2002. \MR{1925796}

\bibitem{Sturmfelsconjecture}
\bysame, \emph{Personal communication}, 2003.

\bibitem{Takeuchi}
M.~Takeuchi, \emph{Free {H}opf algebras generated by coalgebras}, J. Math. Soc.
  Japan \textbf{23} (1971), 561--582. \MR{292876}

\bibitem{Tuttedichromatic}
W.~T. Tutte, \emph{On dichromatic polynominals}, J. Combinatorial Theory
  \textbf{2} (1967), 301--320. \MR{0223272 (36 \#6320)}

\bibitem{Vakil}
R.~Vakil, \emph{Murphy's law in algebraic geometry: badly-behaved deformation
  spaces}, Invent. Math. \textbf{164} (2006), no.~3, 569--590. \MR{2227692}

\bibitem{Weil}
A.~Weil, \emph{Numbers of solutions of equations in finite fields}, Bull. Amer.
  Math. Soc. \textbf{55} (1949), 497--508. \MR{29393}

\bibitem{Welshcomplexity}
D.~J.~A. Welsh, \emph{Complexity: knots, colourings and counting}, London
  Mathematical Society Lecture Note Series, vol. 186, Cambridge University
  Press, Cambridge, 1993. \MR{1245272 (94m:57027)}

\bibitem{White80}
N.~L. White, \emph{A unique exchange property for bases}, Linear Algebra Appl.
  \textbf{31} (1980), 81--91. \MR{570381}

\bibitem{Whitney}
H.~Whitney, \emph{On the {A}bstract {P}roperties of {L}inear {D}ependence},
  Amer. J. Math. \textbf{57} (1935), no.~3, 509--533. \MR{1507091}

\bibitem{WilliamsICM}
L.~Williams, \emph{The positive {G}rassmannian, the amplituhedron, and cluster
  algebras}, Proceedings of the {I}nternational {C}ongress of {M}athematicians
  ({S}t. {P}etersburg), 2022.

\bibitem{Zaslavsky}
T.~Zaslavsky, \emph{Facing up to arrangements: face-count formulas for
  partitions of space by hyperplanes}, Mem. Amer. Math. Soc. \textbf{1} (1975),
  no.~issue 1, 154, vii+102. \MR{0357135 (50 \#9603)}

\end{thebibliography}


\end{document}